\documentclass[final,hidelinks,onefignum,onetabnum,nohypdvips]{siamart220329}

\usepackage{amsfonts}
\usepackage{graphicx,epstopdf} 
\usepackage[caption=false]{subfig} 

\usepackage{mathtools}
\usepackage{bm}
\usepackage{bbm}

\usepackage{algorithmic} 
\Crefname{ALC@unique}{Line}{Lines} 

\usepackage{tikz}
\usetikzlibrary{cd, graphs, positioning,shapes.geometric}
\usetikzlibrary{graphs,decorations.pathmorphing,decorations.markings}

\usepackage{amssymb}
\usepackage{booktabs}
\usepackage{multirow}

\ifpdf
  \DeclareGraphicsExtensions{.eps,.pdf,.png,.jpg}
\else
  \DeclareGraphicsExtensions{.eps}
\fi

\newcommand{\eq}[1]{\begin{align}#1\end{align}}

\definecolor{gold}{rgb}{0,0,1}

\newsiamremark{remark}{Remark}

\title{R-adaptive DeepONet: Learning Solution Operators for PDEs with Discontinuous Solutions Using an R-adaptive Strategy\thanks{Submitted to the editors DATE. \funding{This work was 
			partially supported by the NSF of China grant  12171237, and by the
				Ministry of Science and Technology of China grant 2020YFA0713800, 
				and by the NSFC Major Research Plan - Interpretable and General Purpose Next-generation Artificial Intelligence (Nos. 92270001 and 92370205).}}}

\author{Yameng Zhu\thanks{School of Mathematics,
			Nanjing University, Nanjing 210093, People's Republic of China
  (\email{dg21210022@smail.nju.edu.cn}, \email{wbdeng@nju.edu.cn}).}
\and Jingrun Chen\thanks{School of Mathematical Sciences and Suzhou Institute for Advanced Research, University of Science and Technology of China, Suzhou 215127, People's Republic of China 
  (\email{jingrunchen@ustc.edu.cn}).}
\and Weibing Deng\footnotemark[2]}

\headers{R-adaptive DeepONet}{Yameng Zhu, Jingrun Chen, and Weibing Deng}

\begin{document}
\maketitle
\begin{abstract}
DeepONet has recently been proposed as a representative framework for learning nonlinear mappings between function spaces. However, when it comes to approximating solution operators of partial differential equations (PDEs) with discontinuous solutions, DeepONet poses a foundational approximation lower bound due to its linear reconstruction property. Inspired by the moving mesh (R-adaptive) method, we propose an R-adaptive DeepONet method, which contains the following components: (1) the output data representation is transformed from the physical domain to the computational domain using the equidistribution principle; (2) the maps from input parameters to the solution and the coordinate transformation function over the computational domain are learned using DeepONets separately; (3) the solution over the physical domain is obtained via post-processing methods such as the (linear) interpolation method. Additionally, we introduce a solution-dependent weighting strategy in the training process to reduce the final error. We establish an upper bound for the reconstruction error based on piecewise linear interpolation and show that the introduced R-adaptive DeepONet can reduce this bound. Moreover, for two prototypical PDEs with sharp gradients or discontinuities, we prove that the approximation error decays at a superlinear rate with respect to the trunk basis size, unlike the linear decay observed in vanilla DeepONets. Therefore, the R-adaptive DeepONet overcomes the limitations of DeepONet, and can reduce the approximation error for problems with discontinuous solutions. Numerical experiments on PDEs with discontinuous solutions, including the linear advection equation, the Burgers' equation with low viscosity, and the compressible Euler equations of gas dynamics, are conducted to verify the advantages of the R-adaptive DeepONet over available variants of DeepONet.
\end{abstract}



\section{Introduction}

Many interesting phenomena in physics and engineering are described by partial differential equations (PDEs) whose solutions contain sharp gradient regions or discontinuities.
The most common types of such PDEs are hyperbolic systems of conservation laws  \cite{MR3468916}, such as Euler equations,  inviscid Burgers' equation, etc. It is well-known that solutions of these PDEs develop finite-time discontinuities such as shock waves, even when the initial and boundary data are smooth. 
Other examples include convection-dominated equations, reaction-diffusion equations, and so on.
It is challenging for traditional numerical methods because resolving these discontinuities, such as shock waves and contact discontinuities, requires petite grid sizes. Moreover, characterizing geometric structures, especially in terms of effectively suppressing numerical oscillations near discontinuous interfaces and maintaining the steepness of transition interfaces, is difficult. Specialized numerical methods such as adaptive finite element methods \cite{MR1960405} and discontinuous Galerkin finite element methods \cite{MR3363720} have been successfully used in this context, but their high computational cost limits their wide use.

At the same time, data-driven approaches are becoming a competitive and viable means for solving these challenging problems. 
Deep neural networks (DNNs) have shown promising potential for solving both forward and inverse problems associated with PDEs \cite{MR4270459}. Numerous researchers have explored methods that utilize DNNs for solving PDEs (see \cite{MR4188517, MR4397620} and references therein).

Machine learning for PDEs primarily focuses on learning solutions by training a mapping from the computational domain to the solution. This process, known as the solution parameterization, encompasses techniques such as the Deep Ritz Method \cite{MR3767958}, Deep Galerkin Method \cite{MR3874585}, and Physics Informed Neural Networks (PINNs) \cite{MR3881695, MR4412280}. These methods utilize DNNs to represent the solution and integrate the PDE information into the loss function. The approximate solution is obtained by minimizing the loss function. 
Note that these approaches are tailored to specific instances of PDEs. Consequently, if the coefficients or initial conditions associated with the PDEs change, the model has to be retrained, resulting in poor generalization ability across different PDEs. 

Along another line, there is ongoing work on parameterizing the solution map using DNNs, referred to as operator learning \cite{392253, MR3784416, Khoo, lu2021learning, MR4290514, lifourier}.
In \cite{392253}, Chen and Chen introduced a novel learning architecture based on neural networks, termed operator networks, and demonstrated that these operator networks possess an astonishing universal approximation property for infinite-dimensional nonlinear operators. 
Recently, the authors of \cite{lu2021learning}  replaced the shallow branch and trunk networks in operator networks with DNNs and proposed the deep operator network (DeepONet).
Since proposed, it has been successfully applied to a variety of problems with differential equations \cite{MR4520725, MR4328706, MR4236022, MR4514688}.
In \cite{lifourier}, Li et al. proposed Fourier neural operators based on a nonlinear generalization of the kernel integral representation for some operators and makes use of the convolutional or Fourier network structure. 


Although DeepONets have demonstrated good performance across diverse applications, some studies have pointed out that DeepONets fail to efficiently approximate solution operators of PDEs with sharp gradients or discontinuities \cite{MR4395087, lanthaler2023curse}.
In \cite{MR4395087}, the authors gave a fundamental lower bound on the approximation error of DeepONets and show that there are fundamental barriers to the expressive power of operator learning methods based on linear reconstruction. This is of particular relevance for problems in which the optimal lower bound exhibits a slow decay in terms of the number of basis functions $n$, due to the slow decay of the eigenvalues of the covariance operator.  
To reduce the approximation error, the resolution must be high enough, i.e., we need a large $n$.  However, this may lead to a dramatic increase in computing costs. Therefore, a method with a small approximation error and moderate computational cost is highly desirable.

Variants of DeepONets have been developed to overcome this limitation. Hadorn \cite{hadorn2022shift} investigated the behavior of DeepONet to understand the challenges in detecting sharp features in the target function when the number of basis $n$ is small. They proposed Shift-DeepONet, which adds two neural networks to shift and scale the input function. Venturi and Casey \cite{MR4508150} analyzed the limitations of DeepONet using singular value decomposition and proposed a flexible DeepONet (flexDeepONet) by adding a pre-net and an additional output in the branch net. Seidman et al. \cite{NEURIPS2022_24f49b2a} introduced a nonlinear manifold decoder (NOMAD) framework, utilizing a neural network that takes the output of the branch net as input along with the query location. Recently, Jae Yong Lee et al. \cite{lee2023hyperdeeponet} proposed a HyperDeepONet, which leverages the expressive power of hypernetworks to learn complex operators with a smaller set of parameters. These methods address the limitations of linear reconstruction by modifying the structure of DeepONet, allowing the trunk basis to incorporate information about the input parameters. 

Traditional numerical methods, such as Finite Difference Method and Finite Element Method (FEM), rely on the linear reconstruction using a linear space of basis functions over a predefined mesh to approximate the solution. For solutions with sharp gradients or discontinuities, a fine mesh is needed to resolve local singularities which may lead to significant computational time and data storage. Therefore, researchers have introduced the moving mesh (R-adaptive) method to adaptively and automatically optimize and adjust mesh configurations based on solution characteristics; see \cite{MR2722625} and references therein. The core concept involves adjusting grid distribution through strategic methods without altering the number of mesh grids and their topological connections. This process ensures grids concentrate in regions where solution variations are pronounced. Consequently, this adaptive approach enhances numerical simulation accuracy without increasing computational costs.

To overcome the limitation of linear reconstruction in DeepONet, in this study, we propose a new framework inspired by the moving mesh method, called R-adaptive DeepONet. It employs different learning strategies while maintaining the vanilla structure in the original DeepONet. {To this end,} we introduce a solution-dependent coordinate transformation from the physical domain to the computational domain. The transformed coordinates are then used as the input to the trunk net, similar to traditional R-adaptive methods. This enables adaptive adjustment of basis functions in DeepONet based on the property of the output solution. {Specifically, we first} transform the representation of the output data from the physical domain to the computational domain using the equidistribution principle. This yields two output datasets: the coordinate transform function and the solution over the computational domain. {Second}, we use two DNN models to learn the maps from the input parameters to the coordinate transform function and the solution over the computational domain separately. We emphasize that while learning the forward coordinate transformation from the physical to the computational domain can ensure the injectivity, it retains the singularity of the original solution which is difficult to learn. Therefore, we propose an alternative approach using inverse coordinate transform learning. Although the inverse coordinate transformation does not guarantee a bijection, the functions over the output domain are smoother, making it easier to learn. Given the choice of the inverse coordinate transform, directly predicting the solution value for a given arbitrary coordinate becomes impractical. Thus, we {finally} recover the solution using post-processing methods such as the (linear) interpolation method. {It is worth mentioning that,} according to the error analysis of the operator composition, we introduce two {novel} solution-related weights to the training process of each component. 

We establish an upper bound for the reconstruction error using piecewise linear interpolation and demonstrate that our proposed R-adaptive DeepONet can reduce this bound. Additionally, we rigorously prove that R-adaptive DeepONet can efficiently approximate the prototypical PDEs with sharp gradients or discontinuities. Specifically, the approximation error decays at a superlinear rate with respect to the trunk basis size, while the vanilla DeepONet exhibits at best the linear decay rate \cite{MR4395087}.

To illustrate the effectiveness of our approach, we compare the performance of several DeepONet models for the linear advection equation, the Burgers' equation with low viscosity, and the compressible Euler equations of gas dynamics. The results consistently demonstrate that our R-adaptive DeepONet outperforms vanilla DeepONet and competes effectively with Shift DeepONet.

The remainder of this paper is structured as follows. In \cref{sec:2}, we give a brief introduction to DeepONet. And then discuss the details of the R-adaptive DeepONet in \cref{sec:3}. In \cref{sec:4}, we show some theoretical results. In \cref{sec:5}, we present some numerical results. Finally, some conclusions and comments are given.

\section{Operator Learning and DeepONet}
\label{sec:2}
\subsection{Problem setting}
\label{sec:2.1}
The goal of operator learning is to learn a mapping from one infinite-dimensional function space to another by using a finite collection of observations of input-output pairs from this mapping. We formalize this problem as follows. Let $\mathcal{X}$ and $\mathcal{Y}$ be two Banach spaces of functions defined on bounded domains $D_{\mathcal{X}} \subset \mathbb{R}^{d_{\mathcal{X}}}, D_{\mathcal{Y}} \subset \mathbb{R}^{d_{\mathcal{Y}}}$ respectively and $\mathcal{G}: \mathcal{X} \to \mathcal{Y}$ be a (typically) non-linear map. Suppose we have observations $\{ a^{(i)}, u^{(i)} \}_{i=1}^N$ where $a^{(i)}\sim \mu$ are i.i.d. samples drawn from some probability measure $\mu$ supported on $\mathcal{X}$ and $u^{(i)} = \mathcal{G}(a^{(i)})$. We aim to build an approximation of $\mathcal{G}$ by constructing a parametric map $\mathcal{G}_{\theta}: \mathcal{X}\to \mathcal{Y}$ with parameters $\theta \in \mathbb{R}^{\text{para}}$ such that $\mathcal{G}_{\theta} \approx \mathcal{G}$.

Sometimes the input function space $\mathcal{X}$ can be parameterized by a finite dimensional vector space $\bar{\mathcal{X}}$. Thus, the original objective operator $\mathcal{G} : \mathcal{X} \to \mathcal{Y}$ can also be {equivalently} expressed as {$\bar{\mathcal{G}} : \bar{\mathcal{X}} \to \mathcal{Y}$.} For example, if we consider the mapping from the initial density, velocity, and pressure $(\rho_0, u_0, p_0)$  to the energy $E$ at some time $T$ in the sod shock tube problem, we can parameterize the initial data by the left and right states $(\rho_L, u_L, p_L), (\rho_R, u_R, p_R)$ and the location of the initial discontinuity $x_0$. In this case, the input function space {is equivalent to} a 7-dimensional vector space. For convenience, we still write $\mathcal{G} : \mathcal{X} \to \mathcal{Y}$ instead of distinguishing between $\mathcal{G}$ and $\bar{\mathcal{G}}$.

We are interested in controlling the error of the approximation of the average for $\mu$. In particular, assuming $\mathcal{G}$ is $\mu$-measurable, we aim to control the $L_{\mu}^2(\mathcal{X}; \mathcal{Y})$ Bochner norm of the approximation as follows:
\eq{\label{Bo-error}
\| \mathcal{G} - \mathcal{G}_{\theta} \|_{L^2_{\mu} (\mathcal{X}; \mathcal{Y})} 
:= \mathbb{E}_{a\sim\mu} \| \mathcal{G}(a) - \mathcal{G}_{\theta}(a) \|_{\mathcal{Y}}^2
= \int_{\mathcal{X}} \| \mathcal{G}(a) - \mathcal{G}_{\theta}(a) \|_{\mathcal{Y}}^2\ \text{d} \mu(a).
}

\subsection{A brief introduction to DeepONet}

DeepONets \cite{lu2021learning} present a specialized deep learning architecture for operator learning that encapsulates the universal approximation theorem for operators \cite{392253}. Here we provide a brief introduction to the effective application of DeepONets for learning operators.

To construct a DeepONet, we first need to encode the input parameter function. In \cite{lu2021learning}, the authors use  a fixed collection of training sensors $\{x_1, x_2, \dots, x_m\} \subset D_{\mathcal{X}}$ to encode the input function $a$ by the point values $\mathcal{E}(a) := \mathcal{E}(a(x_1), a(x_2), \dots, a(x_m))$ in $\mathbb{R}^m$. As we mentioned before, sometimes the input function space $\mathcal{X}$ contains a finite-dimensional parameterization and we can encode $a \in \mathcal{X}$ by this parameterization directly. DeepONet is formulated in terms of two neural networks:
\begin{description}
\item[ (1)]  Branch-net ${\bm \beta}$: it maps the point values $\mathcal{E} (a)$ to coefficients 
$${\bm \beta}(\mathcal{E}(a)) = (\beta_1(\mathcal{E}(a)), \dots, \beta_n(\mathcal{E}(a)),$$
 resulting in a mapping
\begin{equation}
\label{eq:deeponetbeta}
{\bm \beta}: \mathbb{R}^m \to \mathbb{R}^n, \quad \mathcal{E}(a) \mapsto {\bm \beta}(\mathcal{E}(a)).
\end{equation}

\item[ (2)]  Trunk-net ${\bm \tau}(y) = (\tau_1(y), \dots, \tau_n(y))$:  it is used to define a mapping
\begin{equation}
\label{eq:deeponettau}
{\bm \tau}: D_{\mathcal{Y}}\to \mathbb{R}^n, \quad y \mapsto {\bm \tau}(y).
\end{equation}
\end{description}

While the branch net provides the coefficients, the trunk net provides the ``basis'' functions in an expansion of the output function of the form
$$
\mathcal{G}^{\text{DON}}(a)(y) = \sum_{k=1}^{n} \beta_k(a)\tau_k(y),\quad a\in \mathcal{X},\ y\in D_{\mathcal{Y}},
$$
with $\beta_k(a) = \beta_k(\mathcal{E}(a))$. The resulting mapping $\mathcal{G}^{\text{DON}}: \mathcal{X}\to \mathcal{Y}$, $a\mapsto \mathcal{G}^{\text{DON}}(a)$ is referred to as the vanilla DeepONet.

\paragraph{Limitation of DeepONet}
Although DeepONets have been proven to be universal within the class of measurable operators \cite{MR4395087}, a fundamental lower bound on the approximation error has also been identified.

\begin{theorem} [Lanthaler et al. 2022, Thm. 3.4]
\label{thm:lowerbound}
Let $\mathcal{X}$ be a separable Banach space, $\mathcal{Y}$ a separable Hilbert space, and let $\mu$ be a probability measure on $\mathcal{X}$. Let $\mathcal{G}: \mathcal{X} \to \mathcal{Y}$ be a Borel measurable operator with $\mathbb{E}_{a\sim \mu} [\|\mathcal{G}(a)\|_{\mathcal{Y}}^2]<\infty$. 
Then the following lower approximation bound holds for any DeepONet $\mathcal{N}^{\text{DON}}$ with trunk-/branch-net dimension $n$:
\begin{equation}
\label{eq:lowerbound}
\mathbin{E}(\mathcal{N}^{\text{DON}}) := \mathbb{E}_{a\sim \mu} [\|\mathcal{N}^{\text{DON}}(a) - \mathcal{G}(a)\|_{\mathcal{Y}}^2]^{1/2} \geq \mathcal{E}_{opt} := \sqrt{\sum_{j>n} \lambda_j},
\end{equation}
where the optimal error $\mathcal{E}_{opt}$ is written in terms of the eigenvalues $\lambda_1 \geq \lambda_2 \geq \dots$ of the covariance operator $\Gamma_{\mathcal{G}_{\# \mu}}:= \mathbb{E}_{u\sim \mathcal{G}_{\# \mu}}[(u\otimes u)]$ of the push-forward measure $\mathcal{G}_{\# \mu}$.
\end{theorem}

The same lower bound applies to any operator approximation of the form $\mathcal{N} (a) = \sum_{k=1}^n \beta_k(a) \tau_k$, where $\beta_k : \mathcal{X} \to \mathbb{R}$ are arbitrary functionals. This bound, for example, also holds for the PCA-Net architecture discussed in \cite{MR3784416} and \cite{MR4290514}. In \cite{LMHM22_1030}, the authors referred to any operator learning architecture of this form as a method with ``linear reconstruction", since the output function $\mathcal{N}(a)$ is restricted to the linear $n$-dimensional space spanned by the $\tau_1, \dots, \tau_n \in \mathcal{Y}$. 

When the eigenvalues $\lambda_1 \geq \lambda_2 \geq \dots$ of the covariance operator $\Gamma_{\mathcal{G}_{\# \mu}}$ decay slowly, approximation using DeepONet may become inaccurate. For instance, the solution operators of advection PDEs and the Burgers' equation are challenging to approximate accurately when using DeepONet with a small number of basis functions $n$ (see \cite{LMHM22_1030}).

\subsection{Variant Models of DeepONet}

Several variants of DeepONet have been developed to overcome its limitations.
Hadorn \cite{hadorn2022shift} proposed the Shift-DeepONet. The main idea is that a scale net ${\bm A} = (A_k)_{k=1}^n$:
$$
{\bm A}: \mathbb{R}^m \to \mathbb{R}^{ n \times d_{\mathcal{Y}} \times d_{\mathcal{Y}} }, \qquad \mathcal{E}(a) \mapsto {\bm A}(a) := \left( A_1(a), A_2(a), \dots A_n(a) \right),
$$
where $A_k(a)$ is matrix-valued functions, and a shift net ${\bm \gamma} = (\gamma_k)_{k=1}^n$, with
$$
 \gamma_k: \mathbb{R}^m \to  \mathbb{R}^{n \times d_{\mathcal{Y}}},  \qquad \mathcal{E}(a) \mapsto {\bm \gamma}(a) := \left(\gamma_1(a), \gamma_2(a), \dots, \gamma_n(a) \right),
$$
to {scale and shift} the input query position $y$, while retaining the DeepONet branch- and trunk-nets ${\bm \beta}, {\bm \tau}$ defined in \cref{eq:deeponetbeta} and \cref{eq:deeponettau}, respectively. The Shift-DeepONet $\mathcal{N}^{\text{sDON}}$ is an operator of the form
$$
\mathcal{N}^{\text{sDON}} (a)(y) = \sum_{k=1}^n \beta_k(a) \tau_k \left( A_k(a)\cdot y + \gamma_k(a) \right).
$$
This approach incorporates the information of the input parameter function $a$ into the trunk basis, allowing the Shift-DeepONet to overcome the limitations of linear reconstruction.

Similar to the Shift-DeepONet, Venturi \& Casey  proposed the flexible DeepONet (flexDeepONet) \cite{MR4508150}, using the additional network, pre-net, to give the bias between the input layer and the first hidden layer, thus introducing the information of $a$ to the trunk basis.
NOMAD \cite{NEURIPS2022_24f49b2a}, developed by Seidman et al., devised a nonlinear output manifold using a neural network that takes the output of the branch net $\{\beta_i\}_{i=1}^n$ and the query location $y$, to overcome the limitation of vanilla DeepONet. Jae Yong Lee et al. went a step further. They used a hypernetwork to share the information of input $a$ to all parameters of the trunk network and proposed a general model HyperDeepONet \cite{lee2023hyperdeeponet}.

All these methods incorporate information from the input function $a$ into the trunk basis to overcome the limitation of linear reconstruction. In practical performance, they do not differ significantly. {To validate the effectiveness of our proposed method, we use Shift DeepONet as a representative among these models to compare with the R-adaptive DeepONet in this paper.}

\section{Proposed Methodology: R-adaptive DeepONet}
\label{sec:3}

\subsection{R-adaptive DeepONet}
\label{sec:3.2}
Many traditional numerical methods rely on linear reconstruction and encounter similar limitations when facing local singularities. R-adaptive methods, also known as moving mesh methods, effectively alleviate these issues. In R-adaptive computations, the number of basis functions remains fixed, but they dynamically adjust based on the problem characteristics. This adaptation reduces errors without significantly increasing computational costs. In \cref{appendix1}, we provide a brief introduction to the R-adaptive method and its associated equidistribution principle.

Inspired by the R-adaptive method, we propose a new learning strategy based on DeepONet for operator learning of PDEs with local singularity, termed R-adaptive DeepONet.

Formally, given $\mathcal{G}: \mathcal{X} \to \mathcal{Y}$, $a \mapsto u(y)$, we introduce a homeomorphism $\tilde{y} = y(\xi) : D_{\mathcal{Y}}\to D_{\mathcal{Y}}, \xi \mapsto y(\xi)$, which maps the computational domain to the physics domain. This allows us to divide the original operator into two new operators as follows:
$$
\mathcal{T}: a \mapsto \tilde{y}(\xi), \quad\text{and}\quad\tilde{\mathcal{G}}: a \mapsto \tilde{u}(\xi) = u(\tilde{y}(\xi)),
$$
where $\mathcal{T}$ maps $a$ to the coordinate transform function $\tilde{y} = y(\xi)$, and $\tilde{\mathcal{G}}$ defines the map to the solution in the computational domain. The original object operator to be learned can be represented as 
\begin{equation}
\label{eq:operatordecompose}
\mathcal{G}(a)(y) = \tilde{\mathcal{G}}(a) \circ \left(\mathcal{T} (a)\right)^{-1} (y),
\end{equation}
where $ \left(\mathcal{T} (a)\right)^{-1}: y\mapsto \xi(y)$ represents the inverse function of $\tilde{y} = y(\xi)$. Since the objective operator $\mathcal{G}(a)$ can be represented by these two operators, we can use two independent DeepONets to learn these two operators as follows: 
$$
\mathcal{T}_{\theta_T}\approx\mathcal{T}: a \mapsto y(\xi) \quad\text{and}\quad
\tilde{\mathcal{G}}_{\theta_G} \approx \tilde{\mathcal{G}} : a \mapsto \tilde{u}(\xi),
$$
where $\theta_T$ and $\theta_G$ represent the parameters of the two models, respectively. For clarity, we will refer to $\mathcal{T}$ as the adaptive coordinate operator and $\tilde{\mathcal{G}}$ as the adaptive solution operator. {The corresponding $\mathcal{T}_{\theta_T}$ and $\tilde{\mathcal{G}}_{\theta_G}$ are termed the adaptive coordinate and adaptive solution DeepONets respectively.} Together, the pair $\{ \mathcal{T}_{\theta_T},  \tilde{\mathcal{G}}_{\theta_G} \}$ is then called an R-adaptive DeepONet system.

Since our approach is data-driven, generating appropriate training data for {the models  $ \mathcal{T}_{\theta_T}$  and $\tilde{\mathcal{G}}_{\theta_G}$} using the equidistribution principle is crucial. Given observations $\{ a^{(i)}, u^{(i)} \}_{i=1}^N$,  we first preprocess the sampled data. This involves determining the corresponding coordinate transform function $y^{(i)}(\xi)$ for each target function $u^{(i)}(y)$ and obtaining the solution on the computational domain $\tilde{u}^{(i)}(\xi)$. As a result, we generate training data sets $\{ a^{(i)}, \{y^{(i)}(\xi_j)\} \}_{i=1}^N$ and $\{ a^{(i)}, \{\tilde{u}^{(i)}(\xi_j)\} \}_{i=1}^N$ for the two independent models, respectively. In this step, we use the mesh generator proposed by Ceniceros and Hou \cite{CENICEROS2001609}. Other mesh generation methods can be found in \cite{MR2722625}.
We emphasize that for problems with discontinuous solutions, the R-adaptive DeepONet needs smaller training datasets than other DeepONets since the output functions $y(\xi)$ and $\tilde{u}(\xi)$ are both smooth, which allow sparser sampling data to capture most features of them.

Here, we choose to learn the mapping from $a$ to $y(\xi)$ instead of $\xi (y)$, since the coordinate transform $y\mapsto \xi(y)$ retains the singularity of the output function, while the inverse $\xi \mapsto y(\xi)$ is relatively smooth, and thus is easier to learn. 

Our goal is to obtain the output function in terms of $y$, but the prediction process yields two functions in terms of $\xi$. To determine the value of $u$ at $y$, we must first find the corresponding $\xi$ and use it as the input of the learned $\tilde{\mathcal{G}}(a)$ to predict the function value.
However, due to the black-box nature of neural networks, deducing the input $\xi$ from output $y$ is challenging. Consequently, post-processing is necessary to make accurate predictions.
After training, we have two independent models mapping the input $a$ to two functions of $\xi$. Given a fixed $a$ and $\xi \in D_{\mathcal{Y}}$, we can get a pair $\{y(\xi), \tilde{u}(\xi)\}$, which forms a mesh grid in the graph of $u = \mathcal{G}(a)$. By using a uniform mesh $\{\xi_j\}$ as the input of the trunk net, we generate a set of points $\{y(\xi_j), \tilde{u}(\xi_j)\}$ that provides a discrete representation of $u$. These points are densely distributed in places where $u$ has singularity, and sparsely distributed in places where $u$ is smooth, hence effectively capturing the function $u$. Using these discrete points, we can reconstruct the output function $u$ by the local interpolation method.

\subsection{Training settings}
\label{sec:3.3}

In \cref{sec:2.1}, we set the target to minimize the $L_{\mu}^2(\mathcal{X}; \mathcal{Y})$ Bochner norm of the approximation~(see \eqref{Bo-error}). In our model, if we assume that the adaptive coordinate operator $\mathcal{T}$ is known and only consider learning the mapping $\tilde{\mathcal{G}}_{\theta_{G}}$, the corresponding approximation error $\mathbin{E}_{\tilde{\mathcal{G}}}$ can be written as follows:
$$
\begin{aligned}
\mathbin{E}_{\tilde{\mathcal{G}}} & = \| \mathcal{G} - \tilde{\mathcal{G}}_{\theta_G} \circ \mathcal{T}^{-1} \|_{L^2_{\mu} (\mathcal{X}; \mathcal{Y})} = \mathbb{E}_{a\sim\mu} \| \mathcal{G}(a) - \tilde{\mathcal{G}}_{\theta_G}(a) (\left(\mathcal{T} (a)\right)^{-1} \|_{\mathcal{Y}}^2\\
&= \mathbb{E}_{a\sim\mu} \int_{D_{\mathcal{Y}}} |\mathcal{G}(a)(y) - \tilde{\mathcal{G}}_{\theta_G}(a)(\left(\mathcal{T} (a)\right)^{-1}(y))|^2 \text{d}y\\
&= \mathbb{E}_{a\sim\mu} \int_{D_{\mathcal{Y}}} |\tilde{\mathcal{G}}(a)(\xi) - \tilde{\mathcal{G}}_{\theta_G}(a)(\xi)|^2 |\text{det} (J(\mathcal{T} (a)(\xi)))| \text{d}\xi.
\end{aligned}
$$
Therefore, in the loss function, we naturally introduce the weight $|\text{det}(J(\mathcal{T} (a)(\xi)))|.$
To prevent this weight from being zero or too large, we modify it to
\begin{equation}
\label{eq:weightg}
w_{\tilde{\mathcal{G}}}(a, \xi) := \min \{ M, \sqrt{1+ |\text{det}(J(\mathcal{T} (a)(\xi)))|^2} \},
\end{equation}
where $M$ is the upper bound we set for this weight, and $|\text{det}(J(\mathcal{T} (a)(\xi)))|$ for weight computing is obtained from the data pre-processing. 
Therefore, in the training process, we aim to minimize the weight empirical loss function
\eq{\label{eq:lossG}
\mathcal{L}_{\tilde{\mathcal{G}}} := \frac{1}{N_1\times N_2} \sum_{k=1}^{N_1} \sum_{j=1}^{N_2} | \tilde{u}_k(\xi_j) - \tilde{\mathcal{G}}_{\theta_G} (a_k)(\xi_j) |^2  w_{\tilde{\mathcal{G}}}(a_k, \xi_j),
}
where $N_1$ denotes the number of sampled inputs $a_k$, and $N_2$ denotes the number of sensors $\xi_j$.
Generally, according to the equidistribution principle, $w_{\tilde{\mathcal{G}}}(a, \xi)$ is relatively small in places where $u$ has singularities. This weighting ensures that the model training is more concentrated over the areas where $u$ is smooth.

In parallel, we can write the approximation error of $\mathcal{T}$ as
$$
\begin{aligned}
\mathbin{E}_{\mathcal{T}} & = \| \mathcal{G} - \tilde{\mathcal{G}} \circ \mathcal{T}_{\theta_T}^{-1} \|_{L^2_{\mu} (\mathcal{X}; \mathcal{Y})} = \mathbb{E}_{a\sim\mu} \int_{D_{\mathcal{Y}}} |\mathcal{G}(a)(y) - \mathcal{G}(a)( \mathcal{T}_{\theta_T}(a) \circ \left(\mathcal{T} (a)\right)^{-1}(y))|^2 \text{d}y\\
& = \mathbb{E}_{a\sim\mu} \int_{D_{\mathcal{Y}}} |\mathcal{G}(a)(y) - \mathcal{G}(a)(\mathcal{T}_{\theta_T}(a)(\xi))|^2 |\text{det} (J(\mathcal{T} (a)(\xi)))| \text{d}\xi\\
& \approx \mathbb{E}_{a\sim\mu} \int_{D_{\mathcal{Y}}} |\nabla \mathcal{G}(a)|^2  |y - \mathcal{T}_{\theta_T}(a)(\xi)|^2 |\text{det} (J(\mathcal{T} (a)(\xi)))| \text{d}\xi.
\end{aligned}
$$
So the corresponding weight can be chosen as 
\begin{equation}
\label{eq:weightt}
w_{\mathcal{T}}(a, \xi) := \min \{ \bar{M}, \sqrt{1+ |\nabla \mathcal{G}(a)|^4|\text{det}(J(\mathcal{T} (a)(\xi)))|^2} \},
\end{equation} 
where $\bar{M}$ is the upper bound we set for this weight.
{The density function is usually of the form $\rho = \sqrt{1 + \beta |\nabla u|^2}$, where $\beta$ is a constant. According to \cref{eq:equidistribution}, we can see that $|\text{det}(J(\mathcal{T} (a)(\xi)))|$ is inversely proportional to $\rho$. So here we can see that $w_{\mathcal{T}}(a, \xi)$ computed according to \cref{eq:weightt} has the opposite performance to $w_{\tilde{\mathcal{G}}}(a, \xi)$. $w_{\mathcal{T}}(a, \xi)$ is small in places where $u$ is smooth and large in places where $u$ has singularities.}

Moreover, for convenience and accuracy of post-processing,  a well-structured mesh is crucial. The coordinate transform functions learned by DeepONet do not inherently guarantee untangling. To prevent mesh tangling, it is essential to ensure that the Jacobian determinant of the transformation function $y(\xi)$ satisfies
$
\text{det}(J(\mathcal{T}_{\theta_T} (a)(\xi))) > 0.
$
Therefore, we incorporate a regularization term into the loss function of the coordinate learning process. The modified loss function becomes
\eq{\label{eq:lossT}\begin{split}
\mathcal{L}_{\mathcal{T}} := &\frac{1}{N_1\times N_2} \sum_{k=1}^{N_1} \sum_{j=1}^{N_2}  \Big[ \lambda_1 | \tilde{y}_k(\xi_j) - \mathcal{T}_{\theta_T} (a_k)(\xi_j) |^2  w_{\mathcal{T}}(a_k, \xi_j)\\
&\qquad\qquad\qquad \qquad+ \lambda_2 \text{ReLU}^2\left(-\text{det}(J(\mathcal{T}_{\theta_T} (a_k)(\xi_j)))\right) \Big],
\end{split}}
where $\lambda_1$ and $\lambda_2$ are regularization parameters, and $w_{\mathcal{T}}(a, \xi)$ represents the weighting factor emphasizing singular regions in $u$.

\section{Theoretical Analysis}
\label{sec:4}
In this section, we provide the theoretical foundation for the effectiveness of our proposed strategy. In traditional numerical methods, R-adaptive methods alleviate the limitations of linear reconstruction by dynamically adjusting the basis functions, thereby reducing approximation errors. Similarly, our proposed R-adaptive DeepONet method can also reduce the errors caused by linear reconstruction. First, we demonstrate the feasibility of the R-adaptive DeepONet in reducing reconstruction errors. Second, we rigorously prove the validity of the proposed method for two prototypical PDEs. In this section, we introduce the shorthand notation $A \lesssim B$ and $B \gtrsim A$ for the inequality $A \leq CB$ and $B \geq CA$, where $C$ denote generic constant independent of the number of trunk net basis functions and the mesh size unless otherwise stated. The notation $A \simeq B$ is equivalent to the statement $A \lesssim B$ and $B \lesssim A$.

\subsection{Reconstruction error of DeepONets}
\subsubsection{Bounds of the reconstruction error}

In \cite{MR4395087}, the authors present a natural decomposition of DeepONets into three components: an encoder $\mathcal{E}$ that maps the infinite-dimensional input space into a finite-dimensional space, an approximator $\mathcal{A}$, often a neural network, maps one finite-dimensional space into another, and a trunk net-induced affine reconstructor $\mathcal{R}$ that maps the finite-dimensional space into the infinite-dimensional output space. 
The total DeepONet approximation error is then decomposed into encoding, approximation, and reconstruction errors.

Suppose $\mathcal{P}_{\tau}: \mathcal{Y} \to {\rm span} \{  \tau_1(y), \dots, \tau_n(y) \}$ is the projection operator that maps the solution function space $\mathcal{Y}$ to the linear span of trunk basis functions. The $L^2$ projection error can be defined as 
$$
\mathbin{E}_{\mathcal{P}_{\tau}} \coloneqq \|\mathcal{P}_{\tau} - {\rm Id} \|_{L^2(\mathcal{G}_{\# \mu})}
= \left( \int_{\mathcal{Y}} \|\mathcal{P}_{\tau} u - u \|^2 {\rm d}(\mathcal{G}_{\#\mu})(u) \right)^{1/2}.
$$
We define the reconstruction error as 
$
\mathbin{E}_{\mathcal{R}} \coloneqq \mathop{\inf}_{\tau} \mathbin{E}_{\mathcal{P}_{\tau}}.
$
The reconstruction error is closely related to the Kolmogorov $n$-width \cite{MR0774404}. According to Kolmogorov $n$-width theory, we can provide a lower bound for the reconstruction error:
$$
\mathbin{E}_{\mathcal{R}} \geq \sqrt{\sum_{j>n} \lambda_j},
$$
where $\lambda_n \geq \lambda_{n+1} \geq \dots$ are defined as in \cref{thm:lowerbound}. This lower bound is fundamental as it reveals that the spectral decay rate for the covariance operator $\Gamma_{\mathcal{G}_{\# \mu}}$ of the push-forward measure essentially determines how low the approximation error of DeepONets can be for a given output dimension $n$ of the trunk nets. 
However in practice, one does not have access to the form of the nonlinear operator $\mathcal{G}$, not to mention the covariance operator $\Gamma_{\mathcal{G}_{\# \mu}}$. Therefore, we aim to derive an upper bound on this error that is easier to be analyzed.

Since the reconstruction error $\mathbin{E}_{\mathcal{R}}$ represents the projection error onto the linear space constructed by the optimal $n$ trunk basis functions in DeepONet, we can compare $\mathbin{E}_{\mathcal{R}}$ with the projection error of a linear reconstruction system constructed using another, possibly non-optimal, set of $n$ basis functions.
A straightforward choice for comparison is the linear finite element reconstruction. In \cite{MR4087799}, the authors proved that a linear finite element function in $\mathbb{R}^d$ with $N$ degrees of freedom can be represented by a ReLU DNN with at most $\mathcal{O}(d)$ hidden layers, and the number of neurons is at most $\mathcal{O}(\kappa^d N)$, where $\kappa \geq 2$ depends on the shape regularity of the underlying finite element grid. 
If we confine the target space to a bounded domain, we can establish an upper bound for the reconstruction error by comparing it to the linear FEM interpolation error.

\begin{lemma}
\label{lem:rerrorleqfemerror}
Suppose that $\Omega \subset D_{\mathcal{Y}}$ is a bounded domain, and $\mathcal{M}_h\subset \Omega$ is a locally convex finite element mesh consisting of a set of simplexes and degrees of freedom $n$. Define the corresponding nodal basis function as $\{  \phi_1(y), \dots, \phi_n(y) \}$, and $\Pi_h$ the interpolation operator on $\mathcal{T}_h$. Then there exists a ReLU-activated trunk net $\bm{\tau}: \mathbb{R}^{d_{\mathcal{Y}}} \to \mathbb{R}^{n}$, with 
$$
\text{depth}(\bm{\tau})  = \mathcal{O}(d_{\mathcal{Y}}),\quad
\text{size}(\bm{\tau}) = \mathcal{O}(\kappa^{d_{\mathcal{Y}}} n),
$$
such that the reconstruction error has the upper bound
$
\mathbin{E}_{\mathcal{R}} \leq \mathbin{E}_{\text{FEM}},
$
where $\mathbin{E}_{\text{FEM}} := \|\Pi_h - {\rm Id} \|_{L^2(\mathcal{G}_{\# \mu})}$ denotes the FEM interpolation error.
\end{lemma}

For the FEM interpolation error, we have the following classical result:
\begin{lemma}
\label{lem:FEMinterperror}
Suppose $\Omega \in \mathbb{R}^d$ is a convex polyhedral domain, and $\mathcal{M}_h\subset \Omega$  is a regular finite element mesh, then
$$
\| u - \Pi_h u\|_{L^2(\Omega)}\ \lesssim\ h^{-2} |u|_{H^2(\Omega)}\quad \forall u\in H^2(\Omega),
$$
where $h = \max_{K\in \mathcal{M}_h} \text{diam}(K)$.
\end{lemma}

Now combining the results of \cref{lem:rerrorleqfemerror}  and \cref{lem:FEMinterperror} yields the following upper bound of the reconstruction error.
\begin{theorem}
\label{thm:errorupperbound}
Suppose $\Omega$ is a bounded convex domain in $\mathbb{R}^{d}$ and $\mathcal{G}$ defines a mapping $\mathcal{G}: \mathcal{X} \to H^2(\Omega)$. Then there exists a trunk net $\bm{\tau}: \mathbb{R}^d \to \mathbb{R}^{n}$, with 
$$
\text{depth}(\bm{\tau})  = \mathcal{O}(d),\quad \text{size}(\bm{\tau}) = \mathcal{O}(\kappa^{d} n),
$$
where $\kappa$ is a constant depending only on $\Omega$, and the associated reconstruction error satisfies
$$
\mathbin{E}_{\mathcal{R}}\ \lesssim\ n^{-\frac{2}{d}} \left( \int_{\mathcal{Y}}|u|^2_{H^2(\Omega)} {\rm d}(\mathcal{G}_{\#\mu})(u) \right)^{1/2}=:\overline{\mathbin{E}_{\mathcal{R}}}.
$$
\end{theorem}

\begin{remark}
Here we obtain an upper bound of the reconstruction error by comparing it with the linear finite element interpolation error on a uniform mesh. However, in general, a uniform mesh is not ideal for finite element interpolation, especially when the objective function has local singularities. In such cases, a very fine mesh is often needed to ensure accuracy, resulting in a large value of $n$, which may lead to an excessively large computational cost.
\end{remark}

We denote the upper bound of the reconstruction error in \cref{thm:errorupperbound} as
$
\overline{\mathbin{E}_{\mathcal{R}}}.
$
In the next part, we will show that our proposed R-adaptive DeepONet has the property of reducing this upper bound of the reconstruction error.

\subsubsection{R-adaptive to lessen the upper bound of the reconstruction error}

Since our proposed R-adaptive DeepONet framework differs from the vanilla DeepONet, we first need to define its reconstruction error.

Suppose $\mathcal{P}_{\tau_{\mathcal{T}}}: \mathcal{Y} \to {\rm span} \{  \tau_{\mathcal{T}, 1}(\xi), \dots, \tau_{\mathcal{T}, n}(\xi) \}$ is the projection operator that maps $\mathcal{Y}$ to the linear span of trunk basis functions of the adaptive coordinate DeepONet $\mathcal{T}_{\theta_T}$. Similarly $\mathcal{P}_{\tau_{\tilde{\mathcal{G}}}}$ is the projection operator that maps $\mathcal{Y}$ to the linear span of trunk basis functions of the adaptive solution DeepONet $\tilde{\mathcal{G}}_{\theta_G}$. Then we define the reconstruction error of the R-adaptive DeepONet as
\begin{equation}
\label{eq:rareconstructionerror}
\mathbin{E}_{\mathcal{R}}^{\text{RA}} \coloneqq \mathop{\inf}_{\tau_{\mathcal{T}}, \tau_{\tilde{\mathcal{G}}}}
\|\mathcal{P}_{\tau_{\tilde{\mathcal{G}}}} \tilde{\mathcal{G}}(a) \circ \left( \mathcal{P}_{\tau_{\mathcal{T}}} \mathcal{T}(a) \right)^{-1} - \text{Id} \|_{L^2(\mu)},
\end{equation}
where $\mathcal{T}, \tilde{\mathcal{G}}$ are introduced in \cref{sec:3.2}, representing the ground adaptive coordinate and solution operators respectively. Moreover, $\mathcal{T}, \tilde{\mathcal{G}}$ satisfy $\tilde{\mathcal{G}}(a) \circ \left(\mathcal{T} (a)\right)^{-1} = \mathcal{G}(a)$. We will explain the reason behind this definition below.
The encoding error arises from the discretization of the input parameter $a$. The approximation error can be viewed as the error associated with learning the linear reconstruction coefficients $\beta(a)$, which is primarily influenced by the branch net. On the other hand, the reconstruction error represents the error due to the inherent linear reconstruction structure in DeepONet, which is affected by the trunk net. Therefore, in our proposed framework with two DeepONets, {we focus solely on the error caused by the trunk net. Combining this error} with \cref{eq:operatordecompose}, we derived the reconstruction errors presented in \cref{eq:rareconstructionerror}.

Assume that $\Omega \subset \mathbb{R}^d$ is polyhedral. $\mathcal{M}_h$ is an affine family of simplicial mesh for $\Omega$, with the reference element $\hat{K}$ being chosen as an equilateral $d$-simplex with unit volume. For any element $K$ in $\mathcal{M}_h$, we denote $F_K: \hat{K} \to K$ as the invertible affine mapping satisfying $K = F_K(\hat{K})$. Then , for any $u\in H^2(\Omega)$, we have
the following error estimate for piecewise linear interpolation:
$$
\|u - \Pi_h u \|_{L^2(\Omega)}^2\ \lesssim\ \sum_{K\in \mathcal{M}_h} \|F'_K\|^4 \cdot |u|_{H^2(K)}^2=:E(u, \mathcal{M}_h),
$$
where $F'_K$ denotes the Jacobian matrix of mapping $F_K$, and  $|u|_{H^2(K)}$ denotes the $H^2$ semi-norm of $u$.
In \cite{MR2722625}, the authors give a lower bound of $E(u,\mathcal{M}_h)$:
$$
E(u,\mathcal{M}_h) \geq N^{-\frac{4}{d}}\left( \sum_{K\in \mathcal{M}_h} |K| \left< u \right>^{\frac{2d}{d+4}}_{H^2(K)} \right)^{\frac{d+4}{d}},
$$
where $N$ is the number of the elements, and $\left< u \right>_{H^2(K)} = \left( \frac{1}{|K|} |u|_{H^2(K)}^2 \right)^{\frac{1}{2}}$. The lower bound can be attained via an optimal mesh $\mathcal{M}_h^{*}$ which equidistributes the density function $\rho_K = \left< u \right>_{H^2(K)}^{\frac{2d}{d+4}}$.
Note that there is a continuous coordinate transformation $x_u: \xi \mapsto x$ from the original uniform simplicial mesh $\mathcal{M}_h$ to the optimal mesh $\mathcal{M}_h^{*}$, which is linear on each element $K\in \mathcal{M}_h$. Through this transform, $\tilde{u} = u\circ x_u$ becomes a function of $\xi$ and can be well approximated on a uniform mesh of $\xi$. Let $\Pi_{\xi}$ denote the linear finite element interpolation operator on the uniform mesh of $\xi$ and $\xi_u: x \mapsto \xi$ the inverse of $x_u$, then we can bound the interpolation error by
\begin{equation}
\label{eq:rerrorofd}
\| u - \left( \Pi_{\xi} \tilde{u} \right) \circ \xi_u \|_{L^2(\Omega)}\ \lesssim\ n^{-\frac{2}{d}} \left( \sum_{K\in \mathcal{M}_h^{*}} |K| \left< u \right>^{\frac{2d}{d+4}}_{H^2(K)} \right)^{\frac{d+4}{2d}}.
\end{equation}
	Here, we have applied the Euler's formula for polyhedra \cite{MR1367739} to address the relationship between the number of elements $N$ and the degrees of freedom $n$.

Based on the previous analysis, we employ ReLU trunk nets to emulate linear finite element space over a predetermined fixed mesh, and then approximate the target function space through function composition. This allows us to derive an upper bound for the reconstruction error. From \cref{eq:rerrorofd} and \cref{lem:rerrorleqfemerror}, it follows that: 
\begin{theorem}
Suppose $\Omega \subset \mathbb{R}^d$ and $\mathcal{G}$ defines a mapping $\mathcal{G}: \mathcal{X} \to H^2(\Omega)$. Then there exists a coordinate transform DeepONet with trunk net $\bm{\tau}_{\mathcal{T}}: \mathbb{R}^d \to \mathbb{R}^{n}$ and a transformed solution DeepONet with trunk net $\bm{\tau}_{\tilde{\mathcal{G}}}: \mathbb{R} \to \mathbb{R}^{n}$ , with 
$$
\begin{aligned}
&\text{depth}(\bm{\tau}_{\mathcal{T}})  = \mathcal{O}(d), &\text{depth}(\bm{\tau}_{\tilde{\mathcal{G}}})  = \mathcal{O}(d),\\
&\text{size}(\bm{\tau}_{\mathcal{T}}) = \mathcal{O}(\kappa^d n), &\text{size}(\bm{\tau}_{\tilde{\mathcal{G}}}) = \mathcal{O}(\kappa^d n),
\end{aligned}
$$
where $\kappa$ is a constant dependent with $d$, and the associated reconstruction error satisfies
$$
\begin{aligned}
\mathbin{E}_{\mathcal{R}}^{\text{RA}} 
\leq \mathbb{E}_{\mathcal{G}_{\#\mu}} \left[ \mathop{\min}_{\mathcal{M}_{h,u}} \mathbin{E}(u, \mathcal{M}_{h,u}) \right]
=\mathbb{E}_{\mathcal{G}_{\#\mu}} \left[\mathbin{E}(u, \mathcal{M}^*_{h,u}) \right]\\
\leq \mathop{\min}_{\mathcal{M}_h} \mathbb{E}_{\mathcal{G}_{\#\mu}} \left[ \mathbin{E}(u, \mathcal{M}_h) \right]\leq \overline{\mathbin{E}_{\mathcal{R}}}.
\end{aligned}
$$
Here $\mathcal{M}_{h,u}$ means the mesh depends on $u$ and $\mathcal{M}^*_{h,u}$ denotes the optimal mesh for each $u$.
\end{theorem}
By adding a coordinate transform (learned by another DeepONet with the same size of the original model), the upper bound for the reconstruction error in the R-adaptive DeepONet may be smaller than that of the vanilla DeepONet. This theorem implies the reduction of the reconstruction error in the R-adaptive DeepONet compared to the vanilla DeepONet.

\subsection{Approximation properties for concrete examples}
The previous subsection theoretically demonstrates that the proposed framework can reduce the upper bound of the reconstruction error, but we do not directly show its advantages over the vanilla DeepONet. This is because the form of the operator $\mathcal{G}$ varies significantly across different problems, making it challenging to use a unified framework for analysis. In this subsection, we select two prototypical PDEs widely used to analyze numerical methods for transport-dominated PDEs. We rigorously prove that the proposed method efficiently approximates operators stemming from discontinuous solutions of PDEs, whereas vanilla DeepONets fail to do so. The chosen PDEs are the linear advection equation and the nonlinear inviscid Burgers' equation, which are the prototypical examples of hyperbolic conservation laws. Detailed descriptions of the exact operators and corresponding approximation results using both {vanilla and our proposed} reconstruction methods are presented below.

\subsubsection{Linear Advection Equation}

Consider the one-dimensional linear advection equation
\begin{equation} 
\partial_t u  + a\partial_x u = 0, \quad u(\cdot, t=0) = \bar{u}
\label{eq:1dlinearadvection}
\end{equation}
on a $2\pi$-periodic domain $\mathbb{T}$, with constant speed $a\in \mathbb{R}$. The underlying operator is $\mathcal{G}_{\rm adv}: L^2(\mathbb{T}) \to L^2(\mathbb{T}), \bar{u} \mapsto \mathcal{G}_{\rm adv}(\bar{u}) \coloneqq u(\cdot, T)$, obtained by solving the PDE \cref{eq:1dlinearadvection} with initial data $\bar{u}$ up to some final time $t= T$. As input measure $\mu \in {\rm Prob}(\mathcal{X})$, we consider random input functions $\bar{u}\sim \mu$ given by the square (box) wave of height $h$, width $w$ and centered at $\zeta$, 
\begin{equation} 
\bar{u}_{\zeta}(x) = h \mathbbm{1}_{[-\omega/2, \omega/2]}(x-\zeta).
\label{eq:linearadvectioninitialdata}
\end{equation}
In the following we let $h=1, w = \pi$, and $\zeta \in [0, 2\pi]$ be uniformly distributed. 

Following \cite{MR4395087}, we observe that the translation invariance of the problem implies that the Fourier basis is optimal for spanning the output space. Given the discontinuous nature of the underlying functions, the eigenvalues of the covariance operator for the push-forward measure decay linearly at most in $n$. Consequently, the lower bound implies a linear decay of error in terms of the number of trunk net basis functions. We derive the following result:
\begin{theorem}
Let $n\in \mathbb{N}$. For any DeepONet $\mathcal{N}^{\rm DON}$ with $n$ trunk-/branch-net output functions, satisfying ${\rm sup}_{\bar{u}\sim \mu}\|\mathcal{N}^{\rm DON}(\bar{u})\|_{L^{\infty}}\leq M< \infty$, we have the lower bound
$$
{\mathcal{E}(\mathcal{N}^{\rm DON})} :=\mathbb{E}_{\bar{u}\sim \mu}[\| \mathcal{G}_{\rm adv}(\bar{u}) - \mathcal{N}^{\rm DON}(\bar{u}) \|^2_{L^2}]^{1/2}\ \gtrsim n^{-1}. 
$$
Consequently, for a given $\epsilon>0$, to achieve $\mathcal{E}(\mathcal{N}^{\rm DON}) \leq \epsilon$ with DeepONet, we need at least $ n\gtrsim \epsilon^{-1}$ trunk and branch net basis functions.
\end{theorem}

In contrast to the previous DeepONet results, we now present an efficient approximation result for R-adaptive DeepONet.
\begin{theorem} 
\label{thm:linearanalysis}
For any $\epsilon >0$, there exist two DeepONets $\mathcal{T}_{\theta_T}, \tilde{\mathcal{G}}_{\theta_G}$, both with $n$ trunk-/branch-net output functions, and assume that $\mathcal{T}_{\theta_T}: [-\pi, \pi] \to [-\pi, \pi], \xi \mapsto x(\xi)$ is bijective. Then the $L^2$-error of the R-adaptive DeepONet system $\{\mathcal{T}_{\theta_T}, \tilde{\mathcal{G}}_{\theta_G}\}$ satisfies
$$
\mathcal{E} := \mathbb{E}_{\bar{u}\sim \mu}[ \| \mathcal{G}_{\rm adv}(\bar{u}) - \tilde{\mathcal{G}}_{\theta_G} (\bar{u}) \circ (\mathcal{T}_{\theta_T}(\bar{u}))^{-1} \|_{L^2}^2 ]^{1/2}\leq \epsilon
$$
with $n\simeq \epsilon^{-2/3}$.
\end{theorem}

The detailed proof, presented in \cref{appendix2}, is based on the fact that the reconstruction error is determined by the approximation error of the optimal reconstruction basis functions. Therefore, if we find a set of basis functions represented by trunk nets that satisfy the error bounds, then the approximation error of the optimal reconstruction basis functions is naturally smaller than the approximation error of this set of basis functions. By construction, we show that the finite element basis functions on a uniform mesh can be represented by trunk nets and satisfies the error bounds. Hence we complete the proof.

\subsubsection{Inviscid Burgers' Equation}
Next, we consider the inviscid Burgers' equation in one-space dimension, which is the prototypical example of nonlinear hyperbolic conservation laws:
\begin{equation}
\partial_t u + \partial_x \left( \frac{1}{2}u^2 \right) = 0, \quad u(\cdot, t=0)= \bar{u},
\label{eq:inviscidburgers}
\end{equation}
on the $2\pi$-periodic domain $\mathbb{T}$. It is well-known that discontinuities in the form of shock waves can appear in finite time even for smooth $\bar{u}$. Consequently, solutions of \cref{eq:inviscidburgers} are interpreted in the sense of distributions and entropy conditions are imposed to ensure uniqueness. Thus, the underlying solution operator is $\mathcal{G}_{\rm Burg} : L^2(\mathbb{T}) \to L^2(\mathbb{T}), \bar{u} \mapsto \mathcal{G}_{\rm Burg}(\bar{u}) = u(\cdot, T)$, with $u$ being the entropy solution of \cref{eq:inviscidburgers} at final time $T$. Given $\zeta \sim \text{Unif}([0, 2\pi])$, we define the random field
$$
\bar{u}_{\zeta}(x) \coloneqq -\sin (x-\zeta),
$$
and we define the input measure $\mu \in {\rm Prob} (L^2(\mathbb{T}))$ as the {law of $\bar{u}_{\zeta}$.} Then, similarly, we can rewrite the underlying operator as $\mathcal{G}_{\rm Burg}: [0, 2\pi] \to L^2(\mathbb{T}), \zeta \mapsto \mathcal{G}_{\rm adv}(\bar{u}_{\zeta}) \coloneqq u_{\zeta}(\cdot, T)$.

Also, translation invariance and local discontinuous can be observed in this problem. Following \cite{MR4395087} we can directly have that
\begin{theorem}
Assume that $\mathcal{G}_{Burg} = u(\cdot, T)$, for $T> \pi$ and $u$ is the entropy solution of \cref{eq:inviscidburgers} with initial data $\bar{u} \sim \mu$. Then the $L^2$-error for any DeepONet $\mathcal{N}^{\rm DON}$ with $n$ trunk-/branch-net output functions is lower-bounded by
$$
{\mathcal{E}(\mathcal{N}^{\rm DON}) :=}\mathbb{E}_{\bar{u}\sim \mu} [\|\mathcal{G}_{\rm Burg}(\bar{u}) - \mathcal{N}^{DON} (\bar{u})\|_{L^2}^2 ]^{1/2}\ \gtrsim n^{-1}.
$$
Consequently, for a given $\epsilon>0$, achieving an error $\mathcal{E}(\mathcal{N}^{\rm DON}) \lesssim \epsilon$ requires at least $n \gtrsim \epsilon^{-1}$.
\end{theorem}

Similar to that in the analysis of linear advection equation, in contrast to the vanilla DeepONet, we have the following result for efficient approximation of $\mathcal{G}_{\rm Burg}$ with R-adaptive DeepONet, whose proof is an almost exact repetition of the proof of \cref{thm:linearanalysis}. 

\begin{theorem} 
Assume that $T>\pi$. For any $\epsilon > 0$, there exist two DeepONets $\mathcal{T}_{\theta_T}, \tilde{\mathcal{G}}_{\theta_G}$, both with $n$ trunk-/branch-net output functions, and assume that $\mathcal{T}_{\theta_T}: [-\pi, \pi] \to [-\pi, \pi], \xi \mapsto x(\xi)$ is bijective. Then the $L^2$-error of the R-adaptive DeepONet system $\{\mathcal{T}_{\theta_T}, \tilde{\mathcal{G}}_{\theta_G}\}$ satisfies
$$
\mathbb{E}_{\bar{u}\sim \mu}[ \| \mathcal{G}_{\rm Burg}(\bar{u}) - \tilde{\mathcal{G}}_{\theta_G} \circ (\mathcal{T}_{\theta_T}(\bar{u}))^{-1} \|_{L^2}^2 ]^{1/2}\leq \epsilon
$$
with $n\simeq \epsilon^{-2/3}$.
\end{theorem}

\section{Numerical Experiments}
\label{sec:5}
In this section, we present several numerical results to evaluate the performance of our proposed R-adaptive DeepONet framework, comparing it with vanilla DeepONets and Shift-DeepONets. We focus on three test problems: Burgers' equation, commonly used to benchmark neural operators; linear advection equations in 1D; and compressible Euler equations in one dimension, which is representative of hyperbolic systems of conservation laws. Through these experiments, we aim to highlight the potential advantages of the proposed framework. The algorithm and definitions of errors in $L^2$ norm are stated in \cref{appendix0}. {In the next part, for simplicity, we sometimes use DON as a shorthand for DeepONet.} \subsection{Linear Advection Equation}
\label{sec:5.1}
\begin{figure}[tbhp] 
\centering
\subfloat[initial condition and final output]{\label{fig:advectiondata(a)}
	\includegraphics[width=0.2\textwidth]{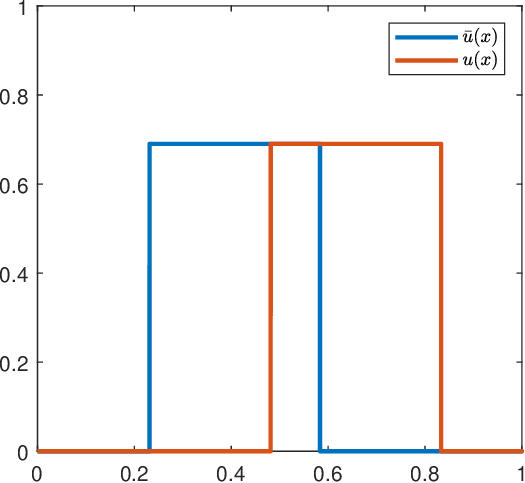}
}
\hspace{2mm}
\subfloat[adaptive coordinate and solution]{\label{fig:advectiondata(b)}
	\includegraphics[width=0.2\textwidth]{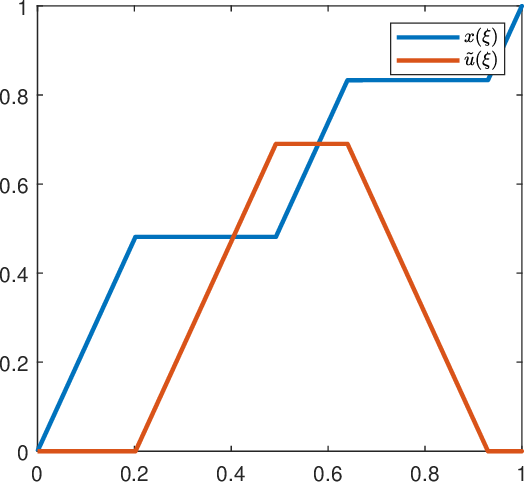}
}
\hspace{2mm}
\subfloat[adaptive coordinate and weight]{\label{fig:advectiondata(c)}
	\includegraphics[width=0.2\textwidth]{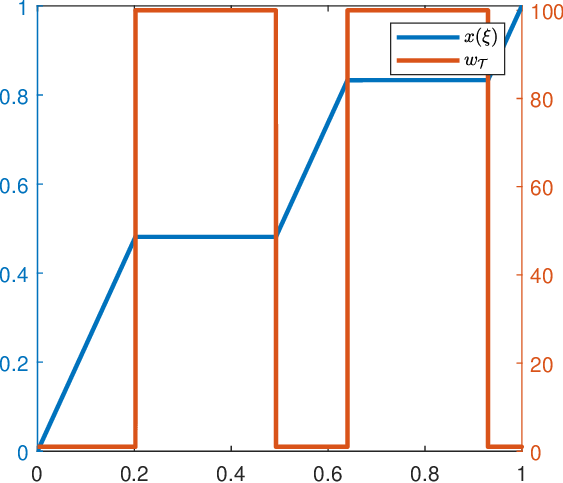}
}
\hspace{2mm}
\subfloat[adaptive solution and weight]{\label{fig:advectiondata(d)}
	\includegraphics[width=0.2\textwidth]{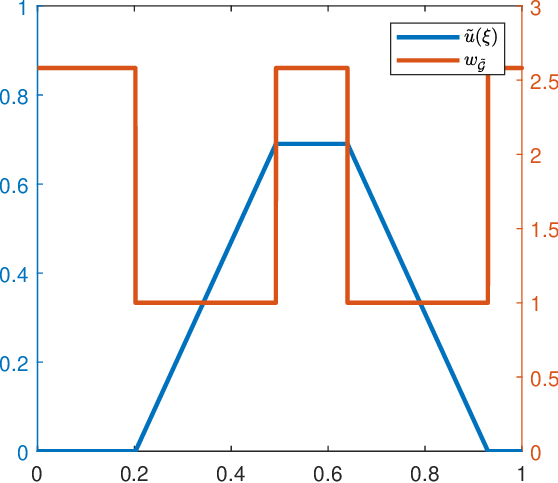}
}
	\caption{Illustration of an example of processed data for advection equation.}
	\label{fig:advectiondata} 
\end{figure}

We take the linear advection equation \cref{eq:1dlinearadvection} as the first example to echo our theoretical analysis in the previous section. Here we set $\Omega =[0,1] $ and $a = 1$. The initial data is given by \cref{eq:linearadvectioninitialdata} corresponding to square waves, with initial heights, widths, and shifts uniformly distributed in $[0.2,0.8]$,  $[0.05, 0.3]$ and $[0, 0.5]$, respectively. We aim to learn the underlying solution operator $\mathcal{G}_{\text{adv}}: \bar{u} \mapsto \mathcal{G}_{\text{adv}}(\bar{u}) = u(\cdot , T=0.25)$, which maps the initial data $\bar{u}$ to the solution at the terminal time $T = 0.25$. Since $\bar{u}$ is controlled by a parameter $\zeta$, the underlying operator is equivalent to $\mathcal{G}: \zeta \mapsto \mathcal{G}_{\rm adv}(\bar{u}_{\zeta}) \coloneqq u_{\zeta}(\cdot, T)$. Therefore, we try to learn the map $\mathcal{G}$ instead of $\mathcal{G}_{\rm adv}$. The training and testing samples of the solutions for vanilla DeepONet and Shift-DeepONet are generated by sampling the underlying exact solution, which are obtained by translating the initial data sampled on $2048$ uniformly distributed grids by $0.25$.  The training data of R-adaptive DeepONet is obtained by preprocessing this batch of data. We use density function 
\begin{equation} 
\label{eq:density1d}
\rho(x) = \sqrt{1+|u'(x)|^2}
\end{equation}
to obtain the equidistributed coordinate transform functions $x(\xi)$ and corresponding adaptive solution functions $\tilde{u}(\xi) = u(x(\xi))$. Moreover, according to \cref{eq:weightg} and \cref{eq:weightt}, we calculate the weights $w_{\tilde{\mathcal{G}}}$ and $w_{\mathcal{T}}$ for the training of the R-adaptive DeepONet. \cref{fig:advectiondata} shows an example of the processed data. As can be seen from Figures {\ref{fig:advectiondata(a)} and \ref{fig:advectiondata(b)}, }the discontinuity in the original data $u(x)$ has been alleviated after preprocessing and has become a smoother transition, and the corresponding coordinate transformation function $x(\xi)$ is also smooth and has no discontinuity. Furthermore, in \cref{fig:advectiondata(c)} and \cref{fig:advectiondata(d)}, we show the calculated weights $w_{\tilde{\mathcal{G}}}$ and $w_{\mathcal{T}}$. We can see that $w_{\tilde{\mathcal{G}}}$ and $w_{\mathcal{T}}$ satisfy certain properties as described in \cref{sec:3}. $w_{\tilde{\mathcal{G}}}$ is relatively small in places where $u$ has singularities, while $w_{\mathcal{T}}$ is just the {opposite.}

To ensure a fair comparison, we used models with similar structures. For vanilla DeepONet,  the main body of Shift-DeepONet, and the two sub-DeepONets in R-adaptive Net, we employed the same architecture: both the branch and trunk nets have 4 layers, each containing 256 neurons. For the scale and shift nets in Shift-DeepONet, we also used a structure of 4 layers with 256 neurons per layer. This approach ensures that the number of parameters in each model remains comparable.

\begin{figure}[tbhp] 
\centering
\subfloat[Validation error during training]{\label{fig:advectionerror(a)}
	\includegraphics[width=0.21\textwidth]{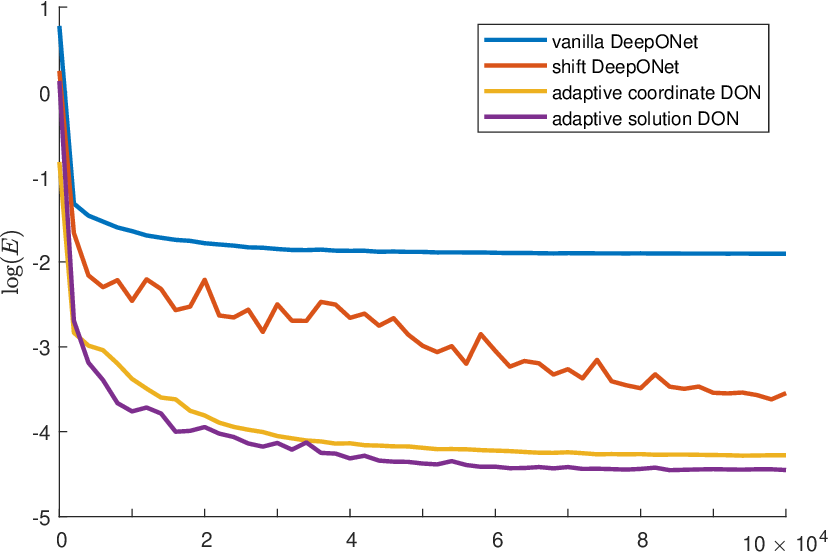}
}
\hspace{2mm}
\subfloat[Eigenvalues of the covariance operators]{\label{fig:advectionsvd(b)}
	\includegraphics[width=0.21\textwidth]{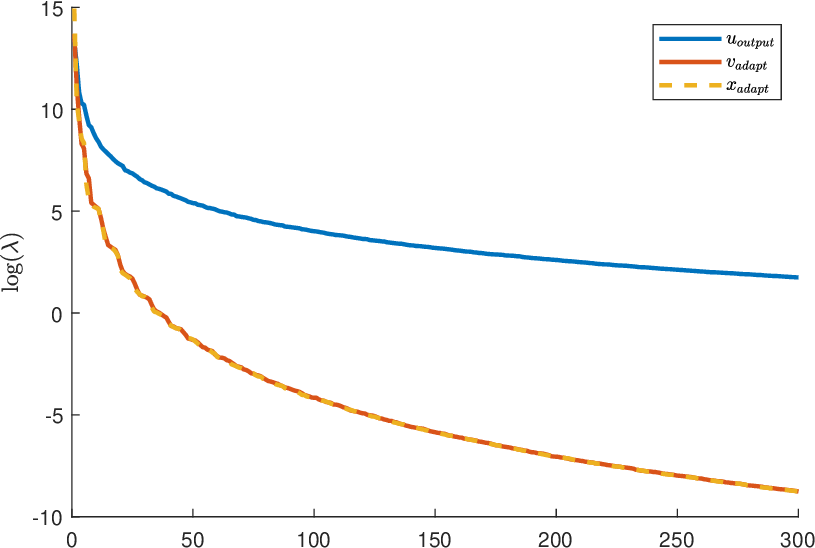}
}
	\caption{Left: the validation error of different models during training; Right: eigenvalues of the covariance operators of different data sets.}
	\label{fig:advectionerror}

\end{figure}

For each model, we use a training set with $1000$ samples and a validation set with $200$ samples. The training is performed with the ADAM optimizer, with learning rate $10^{-3}$ for $100000$ epochs and a learning rate decay of $10\%$. We compute the {relative $L^2$-error} on the validation set every $2000$ epoch. The validation error throughout the training process is shown in \cref{fig:advectionerror(a)}.
It can seen that the validation error of the adaptive solution DON and adaptive coordinate DON decay rapidly, ending up much lower than that of the vanilla DON. Since the validation error indicates the ability of the model to approximate the target dataset, it means that the adaptive solution DON and coordinate DON can approximate their target $\tilde{u}(\xi)$ and $x(\xi)$ well.
The reason for this can be understood by examining the eigenvalues of the covariance matrices of the target datasets. In  \cref{fig:advectionsvd(b)} we show the eigenvalues of the covariance operators of the three data sets, the original solutions $\{u(x)\}$, the processed adaptive solution $\{\tilde{u}(\xi)\}$ and coordinate transform functions $\{x(\xi)\}$. As can be seen from the figure, the eigenvalues of the latter two sets decay much faster than those of the former. From \cref{eq:lowerbound} we know that the corresponding reconstruction error is also smaller. This demonstrates that our preprocessing effectively reduces the lower bound of the reconstruction error, highlighting the feasibility and advantages of our proposed method.

In the testing part, we also use a data set with $200$ samples to calculate the testing error. We use the trained models to predict the solution values at $2048$ grid points uniformly distributed over $[0,1]$ and calculate the approximate $L^2$ error.  For the testing error of R-adaptive DeepONet, we use one-dimensional piecewise-linear interpolation to get the solutions on uniformly distributed grids over $[0,1]$ and then compare it with the exact solutions. 

First, we show the relative testing $L^2$ error of different models trained using output datasets with different sampling densities and verify the advantage of R-adaptive DeepONet that smaller output datasets can be used for good performance.
\cref{tab:advectionerrorresolution} shows the testing errors of models trained using output data sampled on $16$, $32$, $64$ and $128$ uniformly distributed grids over $[0,1]$.
\begin{table}[tbhp]
\footnotesize
\caption{Relative testing $L^2$ error of different models trained using data of different resolutions.}
\label{tab:advectionerrorresolution}
\begin{center}
\begin{tabular}{cccc}
\toprule[2pt] 
Sampling points & vanilla DON & Shift-DON & R-adaptive DON \\
\midrule
16 & $8.17\times 10^{-2}$ & $1.90\times 10^{-2}$ & $6.95\times 10^{-3}$ \\
32 & $4.25\times 10^{-2}$ & $1.68\times 10^{-2}$ & $6.57\times 10^{-3}$ \\
64 & $2.79\times 10^{-2}$ & $1.13\times 10^{-2}$  & $6.96\times 10^{-3}$\\
128 & $2.46\times 10^{-2}$ & $6.37\times 10^{-3}$ & $6.62\times 10^{-3}$\\
\bottomrule[2pt]
\end{tabular}
\end{center}
\end{table}
It is observed that as the number of sampling points increases, the approximation performance of vanilla DeepONet and the Shift-DeepONet improves, resulting in a gradual decrease in testing error. In contrast, R-adaptive DeepONet shows a relative insensitivity to the density of output training data, with the error remaining relatively stable. Therefore, compared to the vanilla DeepONet and the Shift-DeepONet, the accuracy of the R-adaptive DeepONet is less sensistive to the number of sampling points.  
We also note that the R-adaptive DeepONet trained with data sampled on $16$ uniform grids achieves prediction accuracy comparable to that of Shift-DeepONet trained with sampled on $128$ uniform grids. This implies that the proposed method can achieve similar accuracy with a smaller training dataset, hence can reduce the storage requirements during training. This advantage is particularly significant in high-dimensional situations (See the example in \cref{appendix4}).

Next, we present a set of numerical examples to validate the effectiveness of introducing adaptive weights as described in \cref{sec:3.3}. We conducted four groups of experiments using the R-adaptive DeepONet architecture. The training output data is sampled on $2048$ uniformly distributed grids over $[0,1]$. In the first two groups, we 
only learn the adaptive solution operator $\tilde{\mathcal{G}}$
 with upper bounds of the weights: $\bar{w}_{\tilde{\mathcal{G}}} = 1$ and $\bar{w}_{\tilde{\mathcal{G}}} = 2$ respectively. Note that $\bar{w}_{\tilde{\mathcal{G}}} = 1$ indicates training without using weights. In this way, we can demonstrate the effectiveness of introducing weight $w_{\tilde{\mathcal{G}}}$. The results are shown in \cref{tab:advectionerrortable}. It can be seen that the introduction of adaptive weight $w_{\tilde{\mathcal{G}}}$  reduces the approximation error effectively, consistent with the analysis in \cref{sec:3.3}. In the latter two groups of experiments, we change the strategy of learning adaptive coordinates while keeping the part of the adaptive solution unchanged, aiming to show the effectiveness of introducing weight $w_{\mathcal{T}}$. The results also show that the introduction of weight $w_{\mathcal{T}}$ can improve the model's performance.
%

\begin{table}[tbhp]
	\footnotesize
	\caption{Relative testing $L^2$ errors of R-adaptive DeepONet for linear advection equation.}
	\label{tab:advectionerrortable}
	\begin{center}
		\begin{tabular}{ccccc}
			\toprule[2pt] 
			Model &  $\bar{w}_{\tilde{\mathcal{G}}}=1$ \& $x_{\rm ground}$  &$\bar{w}_{\tilde{\mathcal{G}}}=2$ \& $x_{\rm ground}$   & $\bar{w}_{\tilde{\mathcal{G}}}=2$ \& $\bar{w}_{\mathcal{T}}=2$  & $\bar{w}_{\tilde{\mathcal{G}}}=2$ \& $\bar{w}_{\mathcal{T}}=100$  \\
          \midrule
			$Error$  & $5.75\times 10^{-5}$ &
			$3.36\times 10^{-5}$  &
			  $1.69\times 10^{-2}$  &
			  $6.54\times 10^{-3}$  \\
			\bottomrule[2pt]
		\end{tabular}
	\end{center}
\end{table}

\begin{figure}[tbhp] 
\centering
\subfloat[prediction of vanilla DeepONet]{
	\includegraphics[width=0.2\textwidth]{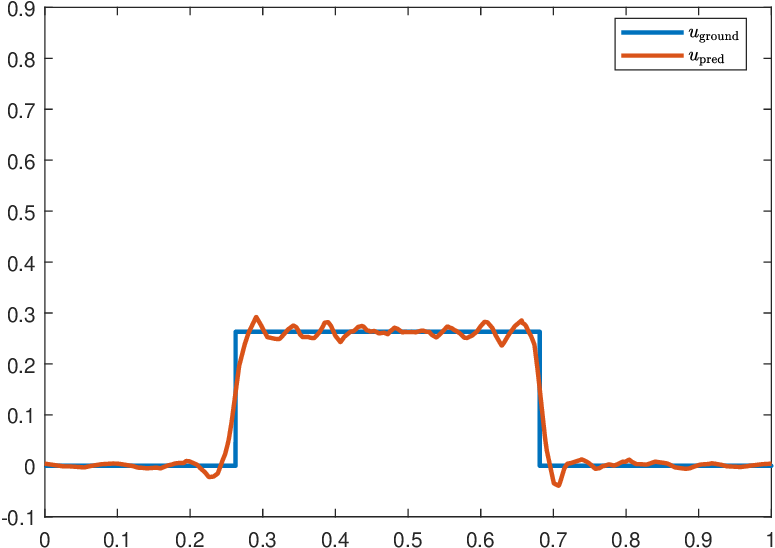}
}
\hspace{2mm}
\subfloat[prediction of Shift-DeepONet]{
	\includegraphics[width=0.2\textwidth]{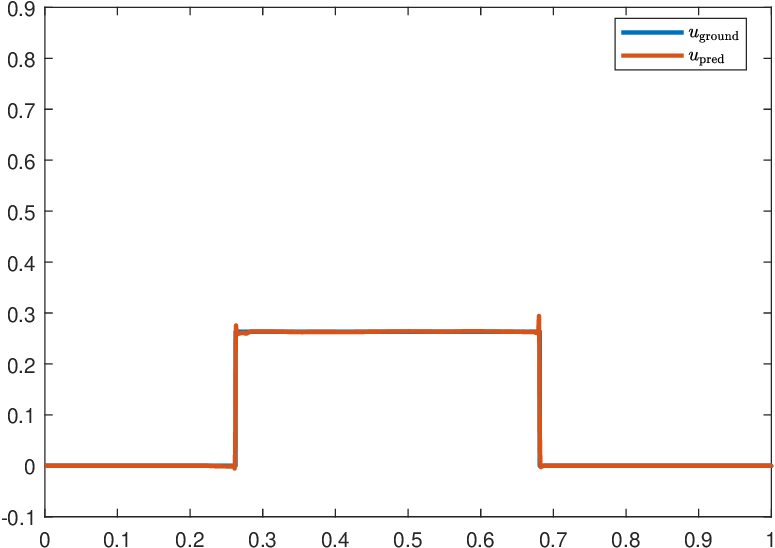}
}
\hspace{2mm}
\subfloat[prediction of R-adaptive DeepONet]{
	\includegraphics[width=0.2\textwidth]{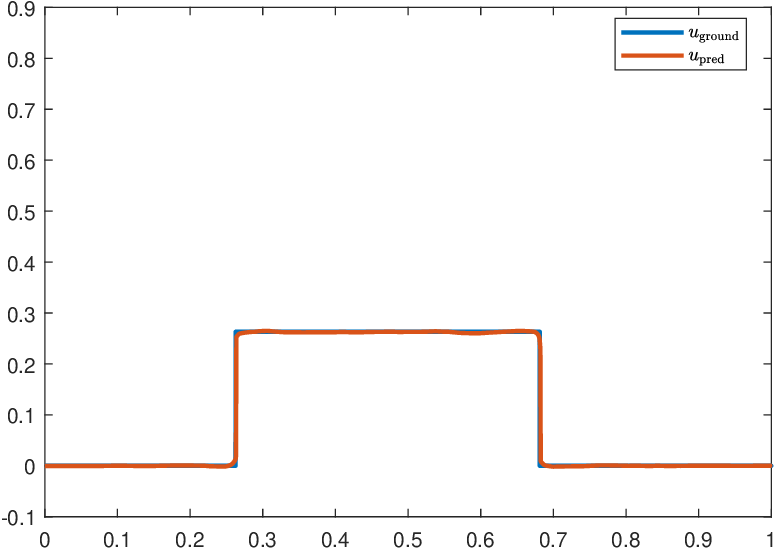}
}

\subfloat[error of vanilla DeepONet]{
	\includegraphics[width=0.2\textwidth]{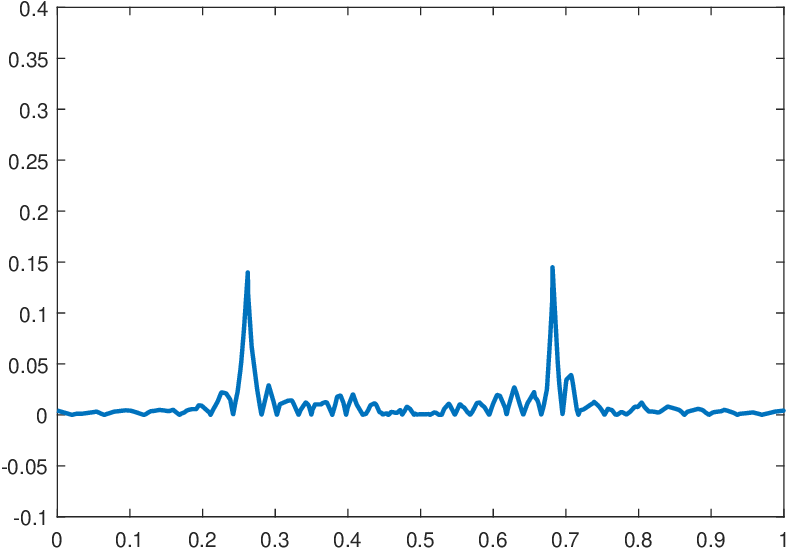}
}
\hspace{2mm}
\subfloat[error of Shift-DeepONet]{
	\includegraphics[width=0.2\textwidth]{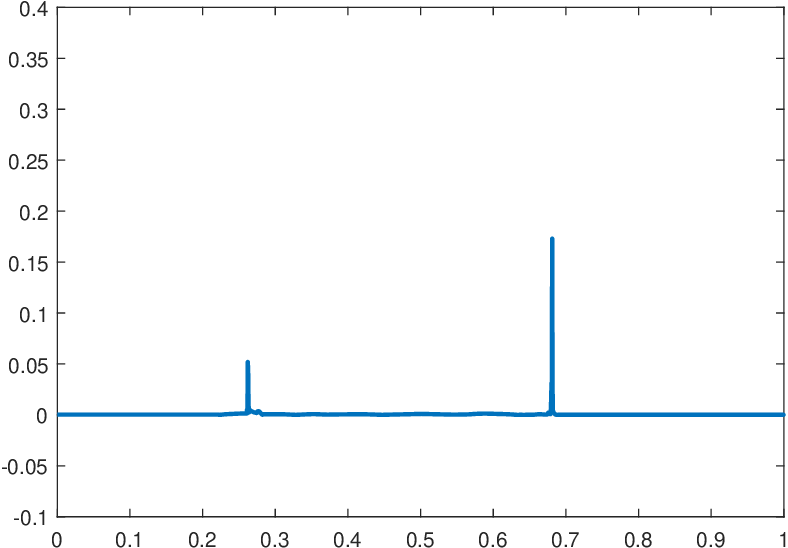}
}
\hspace{2mm}
\subfloat[error of R-adaptive DeepONet]{
	\includegraphics[width=0.2\textwidth]{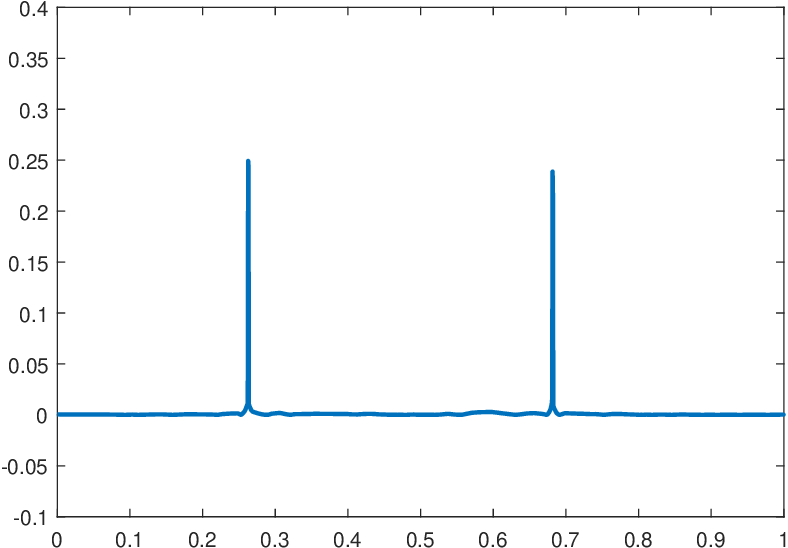}
}
	\caption{An example of the prediction results of the three models for linear advection equation. }
	\label{fig:advectionresult} 
\end{figure}

In the end of this subsection, we provide an example of predictions from different models to visually demonstrate that R-adaptive DeepONet can effectively approximate problems with discontinuities. Here, we use the R-adaptive DeepONet framework with the adaptive solution DeepONet and coordinate DeepONet trained with adaptive weight as shown in the last column in \cref{tab:advectionerrortable}. From \cref{fig:advectionresult}, it can be seen that vanilla DeepONet does not approximate the solution operator well, and its prediction results oscillate wildly due to the existence of discontinuities. Both Shift-DeepONet and R-adaptive DeepONet can grasp the discontinuities and approximate the smooth region well. In addition, Shift-DeepONet leaves small oscillations at the discontinuities, while R-adaptive does not oscillate, but naturally polishes the function. {As a result, both methods still show non-negligible errors at the discontinuities.}

\subsection{Viscous Burgers' Equation}
\label{sec:5.2}
Next, we consider the one-dimensional viscous Burgers' equation
\begin{equation}
\begin{aligned}
\frac{\partial}{\partial t} u(x,t) + \frac{1}{2} \frac{\partial}{\partial x} (u(x,t))^2 
&= \nu \frac{\partial^2}{\partial x^2} u(x,t), &x\in[0,1], t\in[0,1] \\
u(x,0) & = u_0(x), & x\in[0,1]
\end{aligned}
\label{eq:burgersequation}
\end{equation}
with periodic boundary conditions and a fixed viscosity $\nu$.

When the viscosity coefficient is large, the solution of the Burgers' equation will not exhibit significant singularities. However, as the viscosity coefficient decreases, the solution gradually approaches that of the corresponding inviscid Burgers' equation, resulting in regions with large gradients. In this experiment, we use several different viscosity coefficients such as $\nu = 5\times  10^{-2}, 10^{-2}, 10^{-3}$, and $10^{-4}$. Our goal is to learn the solution operator mapping initial conditions $u(x, 0)$ to the solution at $T=1$. 

\begin{figure}[tbhp] 
\centering
\subfloat[initial condition and solution at $t=1$ for $\nu = 10^{-2}$]{
	\includegraphics[width=0.2\textwidth]{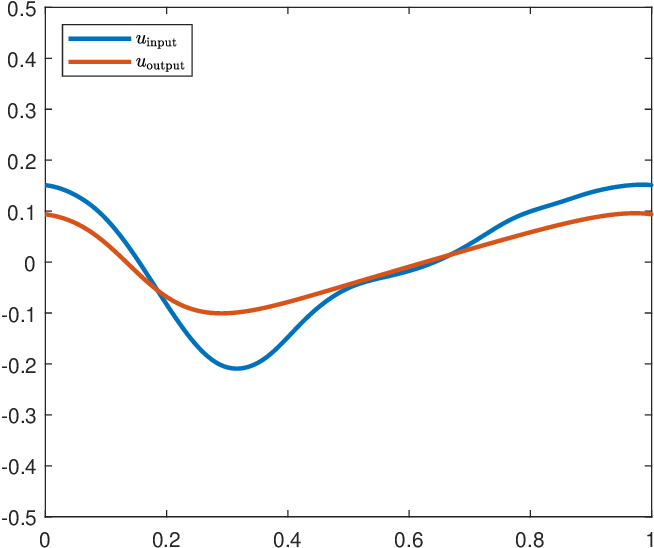}
}
\hspace{2mm}
\subfloat[adaptive coordinate and weight for $\nu = 10^{-2}$]{
	\includegraphics[width=0.21\textwidth]{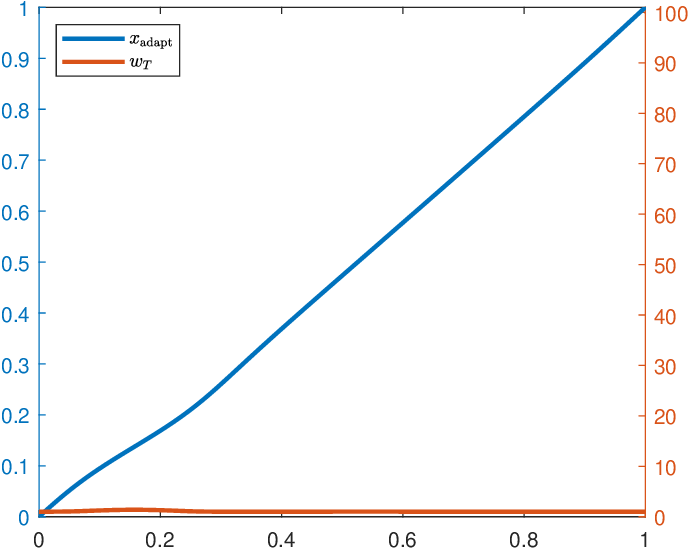}
}
\hspace{2mm}
\subfloat[adaptive solution and weight for $\nu = 10^{-2}$]{
	\includegraphics[width=0.21\textwidth]{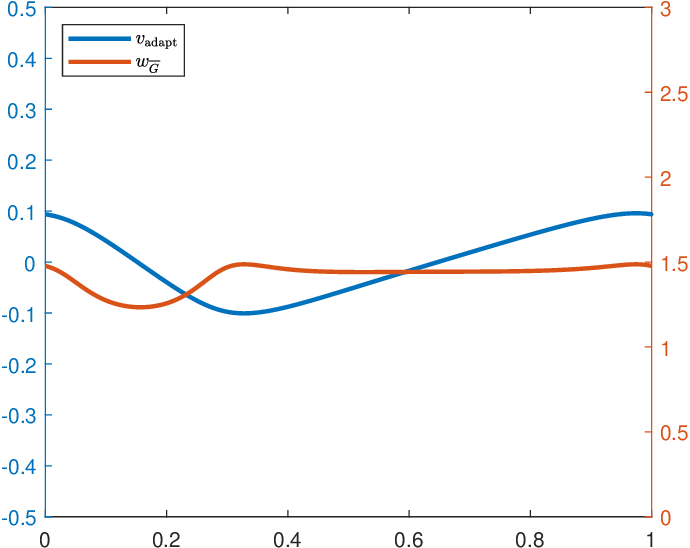}
}

\subfloat[initial condition and solution at $t=1$ for $\nu = 10^{-3}$]{
	\includegraphics[width=0.2\textwidth]{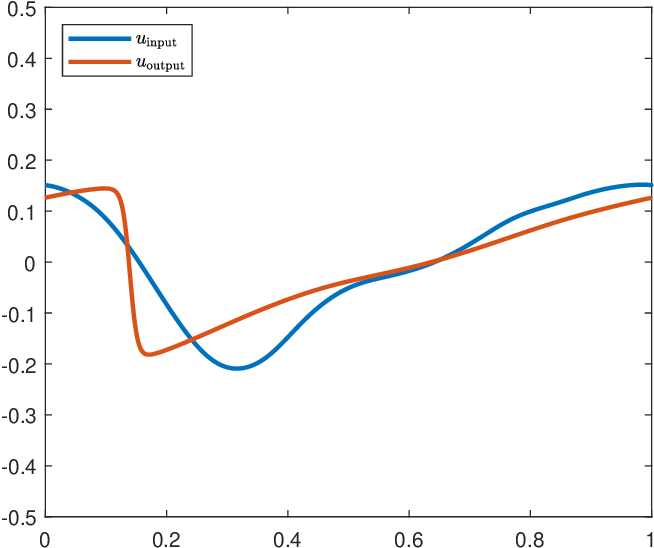}
}
\hspace{2mm}
\subfloat[adaptive coordinate and weight for $\nu = 10^{-3}$]{
	\includegraphics[width=0.21\textwidth]{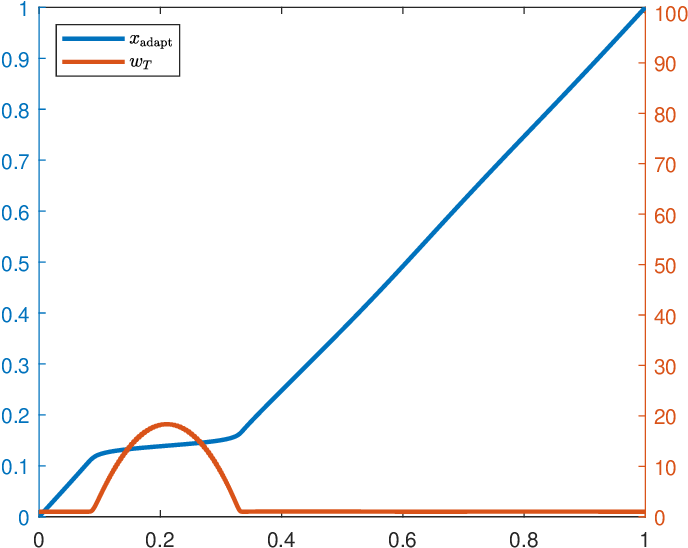}
}
\hspace{2mm}
\subfloat[adaptive solution and weight for $\nu = 10^{-3}$]{
	\includegraphics[width=0.21\textwidth]{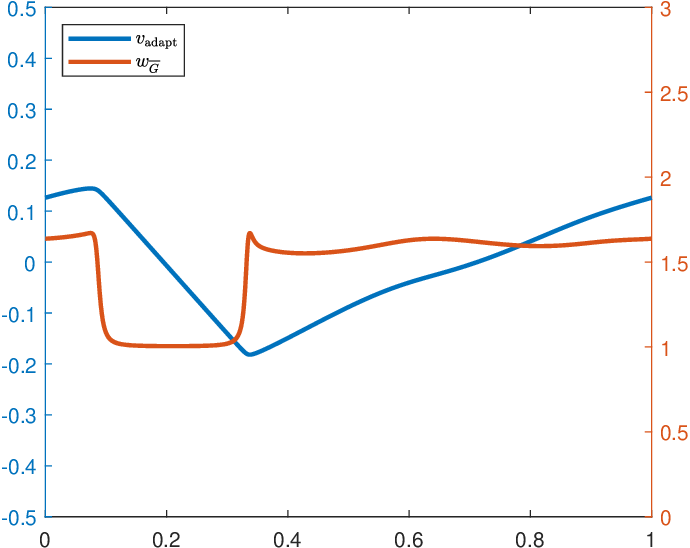}
}

\subfloat[initial condition and solution at $t=1$ for $\nu = 10^{-4}$]{
	\includegraphics[width=0.2\textwidth]{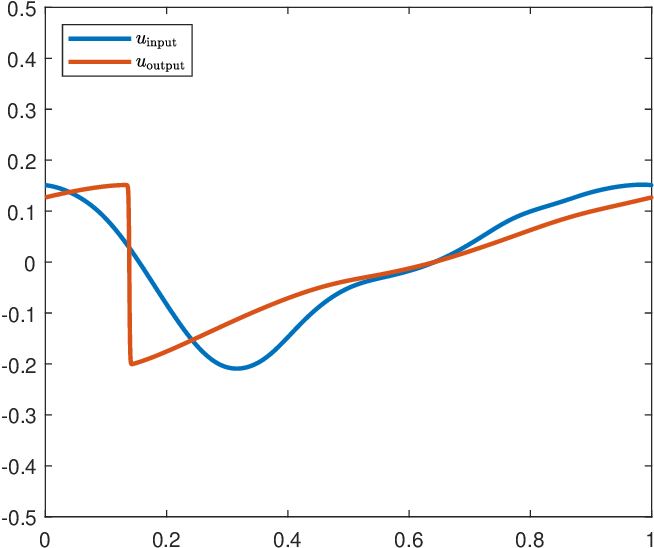}
}
\hspace{2mm}
\subfloat[adaptive coordinate and weight for $\nu = 10^{-4}$]{
	\includegraphics[width=0.21\textwidth]{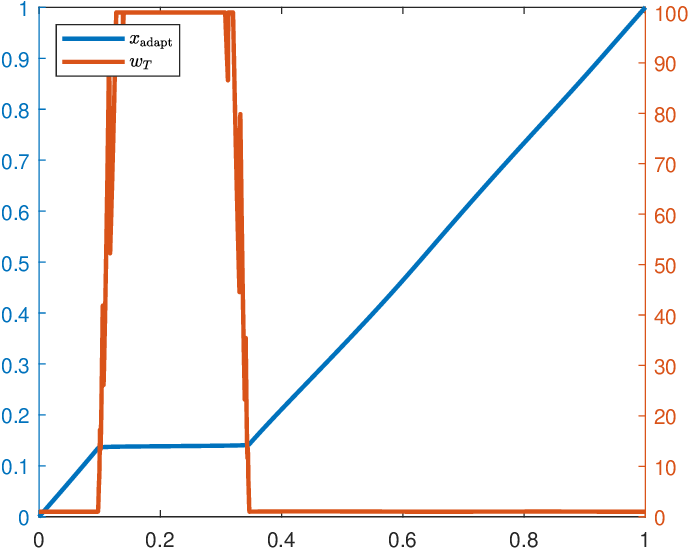}
}
\hspace{2mm}
\subfloat[adaptive solution and weight for $\nu = 10^{-4}$]{
	\includegraphics[width=0.21\textwidth]{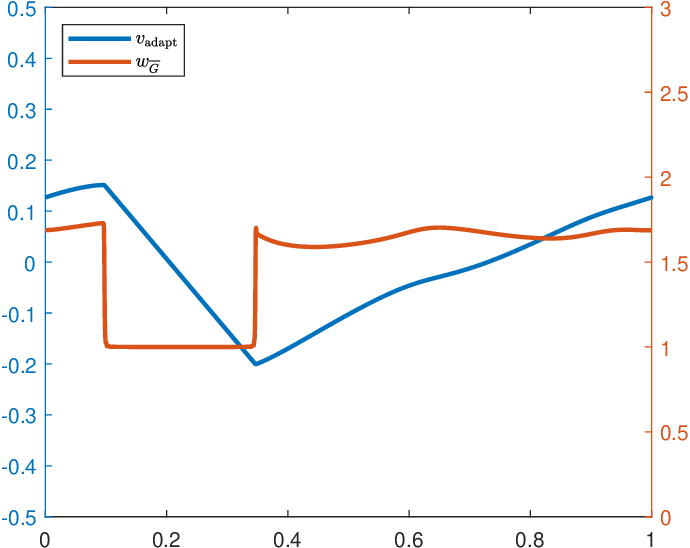}
}
	\caption{Illustration of an example of processed data for Burgers' equation.}
	\label{fig:burgersdata} 
\end{figure}

To obtain a set of training data, we randomly sample 1000 input functions from a Gaussian Random Field (GRF) $~\mathcal{N}(0, 25^2 (-\Delta + 5^2 I)^{-4})$ and solve the Burgers' equation using the Chebfun package with a spectral Fourier discretization and a fourth-order stiff time stepping scheme with a time-step size of $10^{-4}$. we generate test data sets by sampling another 200 input functions from the same GRF. On the input side, we sample the initial data on a uniformly distributed grid of $128$ points over $[0,1]$ as the input parameters for training the models. The data preprocessing is similar to the previous test, and the density function \cref{eq:density1d} is also used. As before, we show some examples of the processed data in \cref{fig:burgersdata}. When the viscosity coefficient is relatively large, the solution does not exhibit singularities. In this case, the adaptive solution obtained through preprocessing is close to the original data, and the coordinate transformation is approximately an identity mapping. 
However, as the viscosity coefficient decreases, the adaptive solution obtained through preprocessing becomes smoother and free of singularities compared to the original data. 
Additionally, we have shown the {graphs} of the adaptive weights, whose properties are consistent with our analysis in \cref{sec:3.3}.

We use the same network structures and training strategies as in the previous experiment. The testing errors for different operator learning strategies are presented in \cref{tab:burgerserrortable}. As shown in the table, vanilla DeepONet approximates the solution operator well when the viscosity coefficient is large.
As the viscosity coefficient decreases and the solution exhibits local singularities, the performance of vanilla DeepONet degrades. In contrast, R-adaptive DeepONet performs better than vanilla DeepONet at low viscosity levels and even achieves smaller relative errors than Shift-DeepONet.
This may indicate that R-adaptive DeepONet has an advantage over Shift-DeepONet in approximating problem whose solution exhibits large gradients rather than discontinuities, such as convection-dominated diffusion equations. We will explore this in future work.

\begin{table}[tbhp]
\footnotesize
\caption{Relative testing $L^2$ error of different models for Burgers' equation.}
\label{tab:burgerserrortable}
\begin{center}
\begin{tabular}{ccccc}
\toprule[2pt] 
{ Model} &{$\nu = 5\times 10^{-2}$} &{$\nu = 10^{-2}$}  & {$\nu = 10^{-3}$} &{$\nu = 10^{-4}$} \\
\midrule 
vanilla DeepONet& $5.41\times 10^{-5}$ & $3.93\times 10^{-4}$ &  $1.00\times 10^{-2}$ & $3.45\times 10^{-2}$  \\
Shift-DeepONet& $1.58\times 10^{-4}$  & $6.32\times 10^{-4}$ &  $1.15\times 10^{-2}$  & $3.93\times 10^{-2}$\\
R-adaptive DON& $1.31\times 10^{-4}$ & $4.90\times 10^{-4}$ &  $8.35\times 10^{-3}$  & $2.44\times 10^{-2}$\\
\bottomrule[2pt] 
\end{tabular}
\end{center}
\end{table}

\begin{figure}[tbhp] 
\centering
\subfloat[vanilla DON for $\nu = 10^{-2}$]{
	\includegraphics[width=0.2\textwidth]{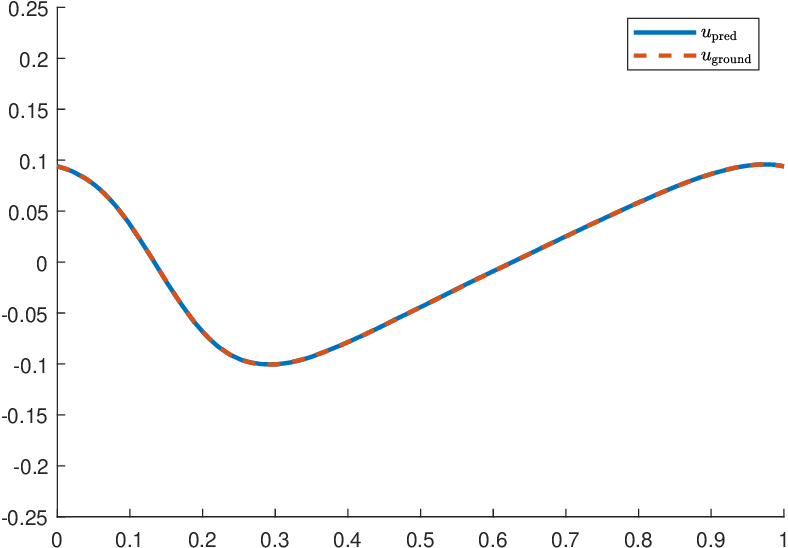}
}
\hspace{2mm}
\subfloat[Shift DON for $\nu = 10^{-2}$]{
	\includegraphics[width=0.2\textwidth]{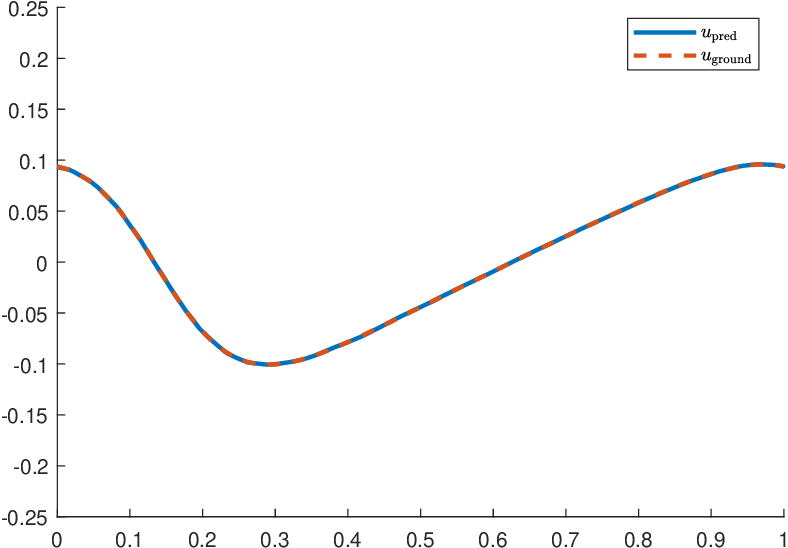}
}
\hspace{2mm}
\subfloat[R-adaptive DON for $\nu = 10^{-2}$]{
	\includegraphics[width=0.2\textwidth]{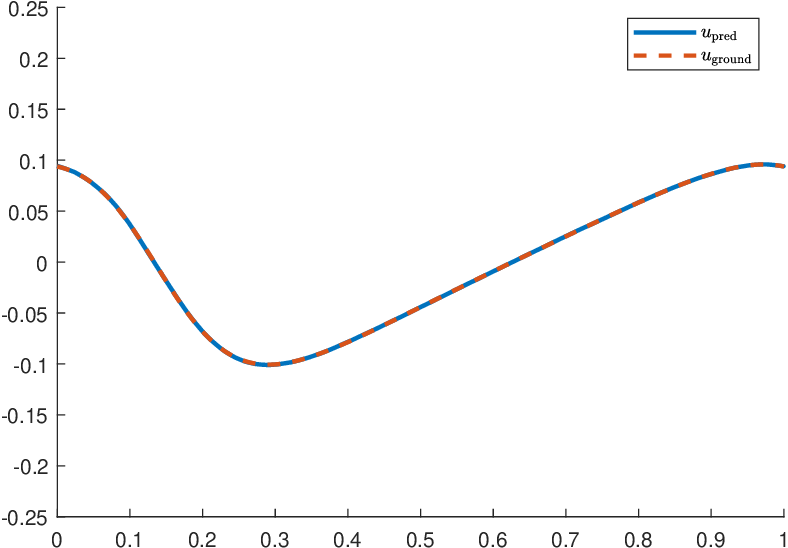}
}

\subfloat[vanilla DON for $\nu = 10^{-3}$]{
	\includegraphics[width=0.2\textwidth]{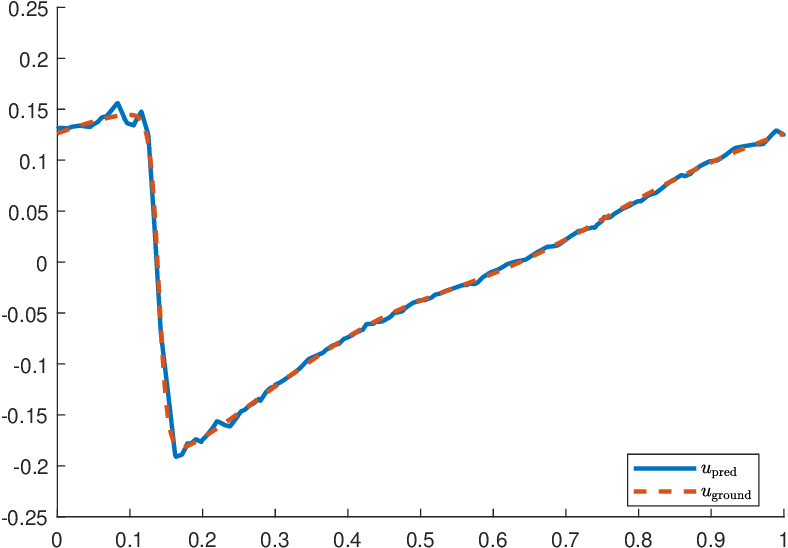}
}
\hspace{2mm}
\subfloat[Shift DON for $\nu = 10^{-3}$]{
	\includegraphics[width=0.2\textwidth]{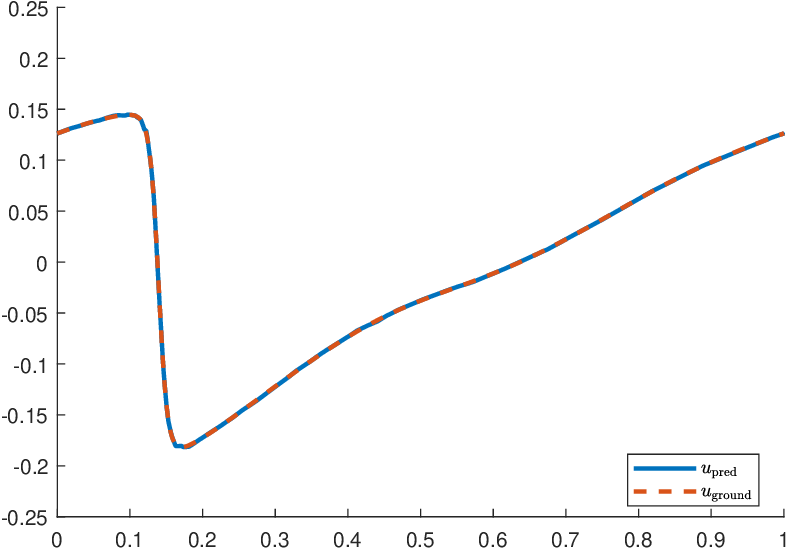}
}
\hspace{2mm}
\subfloat[R-adaptive DON for $\nu = 10^{-3}$]{
	\includegraphics[width=0.2\textwidth]{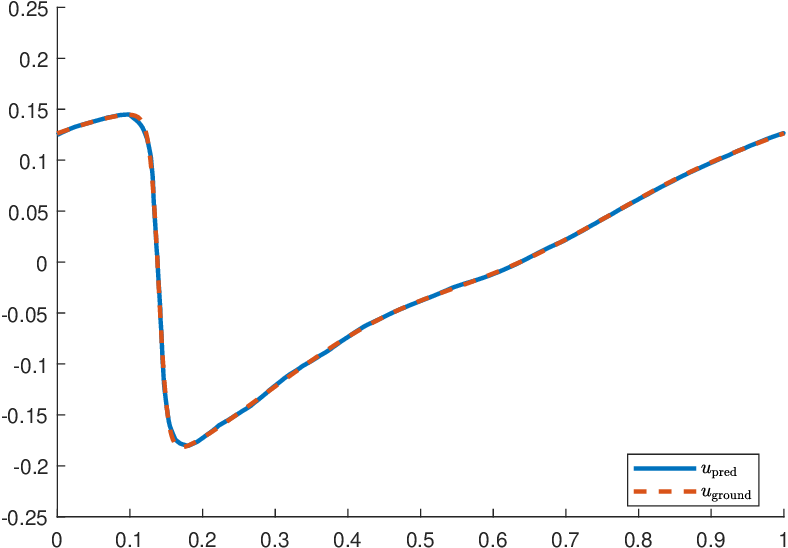}
}

\subfloat[vanilla DON for $\nu = 10^{-4}$]{
	\includegraphics[width=0.2\textwidth]{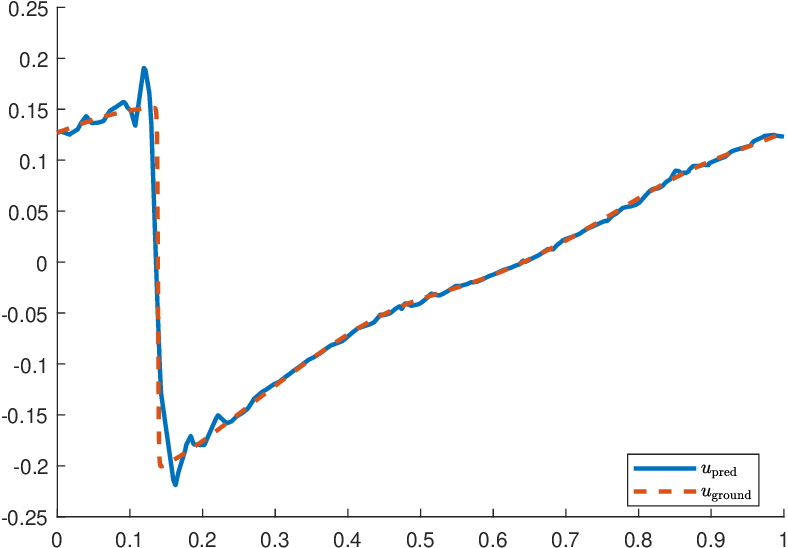}
}
\hspace{2mm}
\subfloat[Shift DON for $\nu = 10^{-4}$]{
	\includegraphics[width=0.2\textwidth]{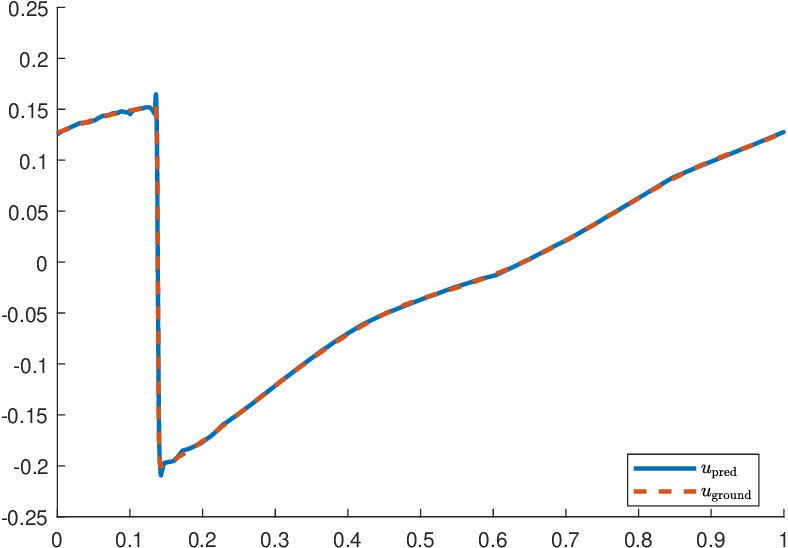}
}
\hspace{2mm}
\subfloat[R-adaptive DON for $\nu = 10^{-4}$]{
	\includegraphics[width=0.2\textwidth]{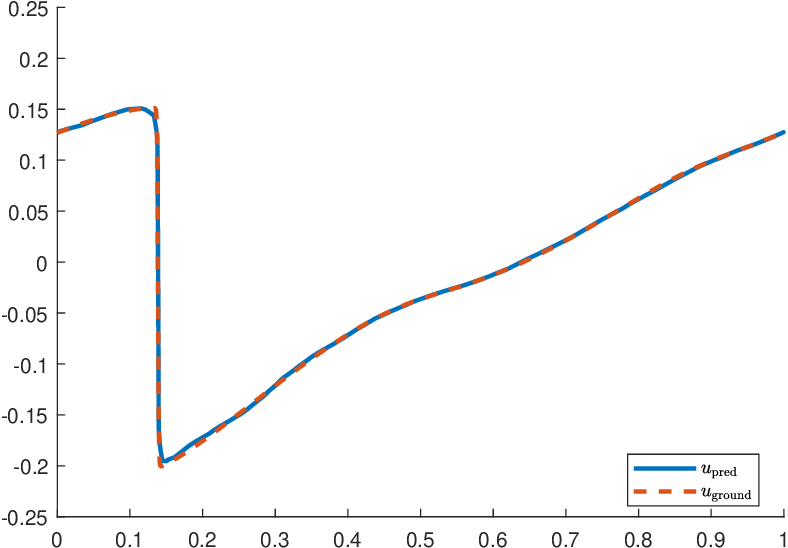}
}

	\caption{Illustration of an example of processed data for Burgers' equation.}
	\label{fig:burgersfigure} 
\end{figure}

To illustrate the prediction results more intuitively, we present some prediction examples in \cref{fig:burgersfigure}. As seen in the figures, vanilla DeepONet performs well in approximating the solution when the viscosity coefficient is large, and both Shift-DeepONet and R-adaptive DeepONet also provide accurate predictions. However, as the viscosity coefficient decreases, the solution of Burgers' equation develops a large local gradient, causing the solution predicted by vanilla DeepONet to oscillate, especially near the singularity region. In contrast, both Shift- and R-adaptive DeepONets capture the local singularity characteristics effectively. Additionally, as noted in \cref{sec:5.1}, Shift-DeepONet produces minor oscillations while R-adaptive DeepONet polishes the solutions around the singularity.

 In addition to predicting one-dimensional functions, our proposed R-adaptive DeepONet is also suitable for higher dimensional cases. We provide an example in \cref{appendix4} to verify this point.

\subsection{Shock Tube}
In this subsection, we consider the motion of an inviscid gas described by the Euler equations of aerodynamics. The governing equations can be written as
$$
\begin{pmatrix}
\rho\\
\rho u\\
E
\end{pmatrix}_t
+ 
\begin{pmatrix}
\rho u\\
\rho u^2 + p\\
(E+p)u
\end{pmatrix}_x = 0,
$$
with $\rho, u$ and $p$ denoting the fluid density, velocity, and pressure. $E$ represents the total energy per unit volume: $
E= \frac{1}{2}\rho u^2 + \frac{p}{\gamma -1},
$
where $\gamma = c_p/c_v$ is the gas constant which equals to $1.4$ for a diatomic gas considered here.

We restrict the equation to $D=[-5,5]$ and consider the initial data corresponding to a shock tube of the form
$$
\rho_0 = \left\{
\begin{aligned}
&\rho_L &x\leq x_0\\
&\rho_R&x> x_0
\end{aligned}
\right.
\quad
u_0 = \left\{
\begin{aligned}
&u_L &x\leq x_0\\
&u_R&x> x_0
\end{aligned}
\right.
\quad
p_0 = \left\{
\begin{aligned}
&p_L &x\leq x_0\\
&p_R&x> x_0
\end{aligned}
\right.
$$
parameterized by the left and right states $(\rho_L, u_L, p_L)$, $(\rho_R, u_R, p_R)$, and the location of the initial discontinuity $x_0$. As proposed in Lye et al. \cite{MR4078984}, these parameters are, in turn, drawn from the measure 
$$
\begin{aligned}
&\rho_L = 0.75 + 0.45g(z_1),\quad u_L = 0.5+0.5g(z_3),\quad  p_L=2.5+ 1.6g(z_4),\\
&\rho_R = 0.4 + 0.3g(z_2),\quad u_R = 0,\quad  p_R= 0.375+ 0.325g(z_5), \quad x_0 = 0.5g(z_6),& & & &
\end{aligned}
$$
with $z = [z_1, \dots, z_6] \sim \text{Unif}([0,1]^6)$ and $g(z) = 2z-1$. We aim to approximate the operator $\mathcal{G}: [\rho_0, \rho_0 u_0, E_0] \mapsto E(1.5)$. As in the previous subsections, we simplify this mapping to $\mathcal{G}: z \mapsto E(1.5)$.

\begin{figure}[tbhp] 
\centering
\subfloat[initial $\rho_0$ and $u_0$]{
	\includegraphics[width=0.2\textwidth]{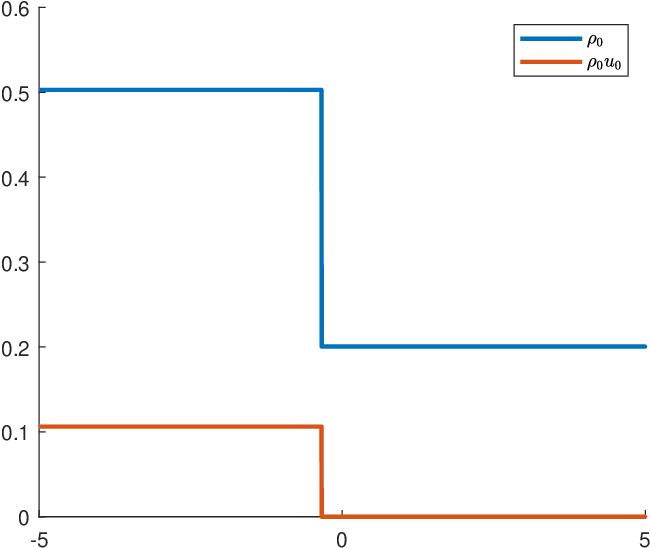}
}
\hspace{2mm}
\subfloat[initial $E_0$ and the final output $E$]{
	\includegraphics[width=0.2\textwidth]{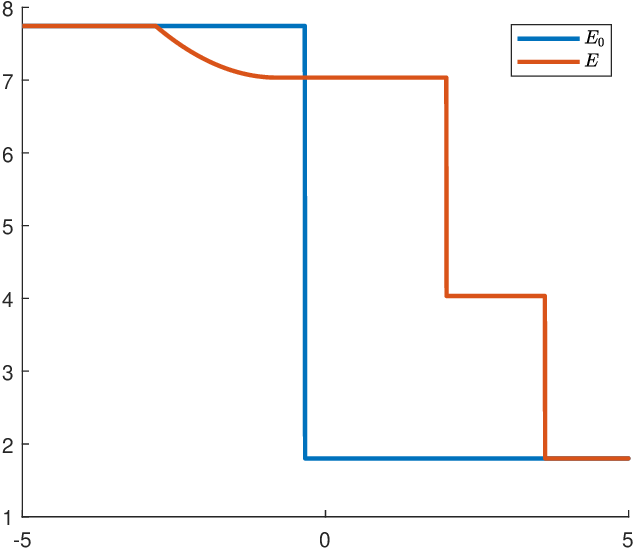}
}
\hspace{2mm}
\subfloat[adaptive coordinate and weight]{
	\includegraphics[width=0.2\textwidth]{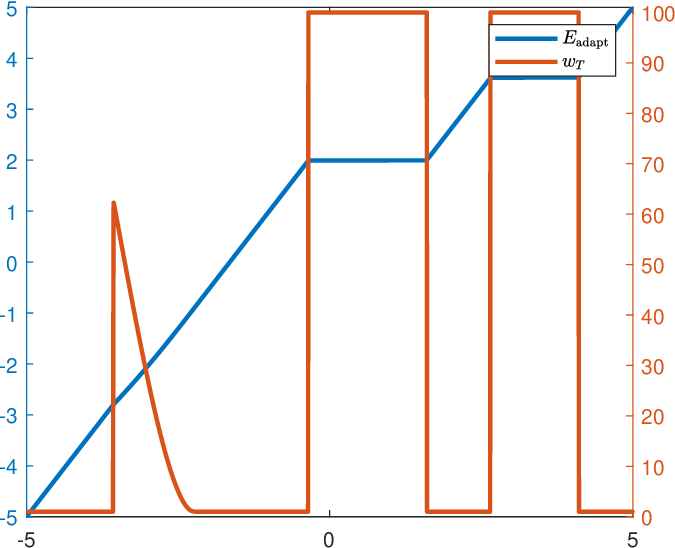}
}
\hspace{2mm}
\subfloat[adaptive solution and weight]{
	\includegraphics[width=0.2\textwidth]{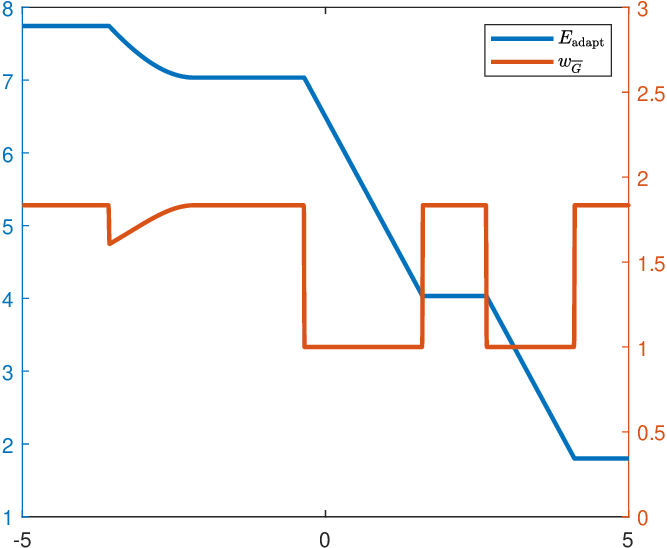}
}
	\caption{Illustration of an example of processed data for Sod shock tube problem.}
	\label{shocktubeproposeddata} 
\end{figure}

The training (and testing) output is generated through the analytic method in \cite{toro1999numerical}. The rest of the concretes are similar to those in subsections \ref{sec:5.1} and \ref{sec:5.2}. As in the previous subsections, we first show an example of the processed data in \cref{shocktubeproposeddata}.

\begin{figure}[tbhp] 
\centering
\subfloat[prediction of vanilla DON]{
	\includegraphics[width=0.25\textwidth]{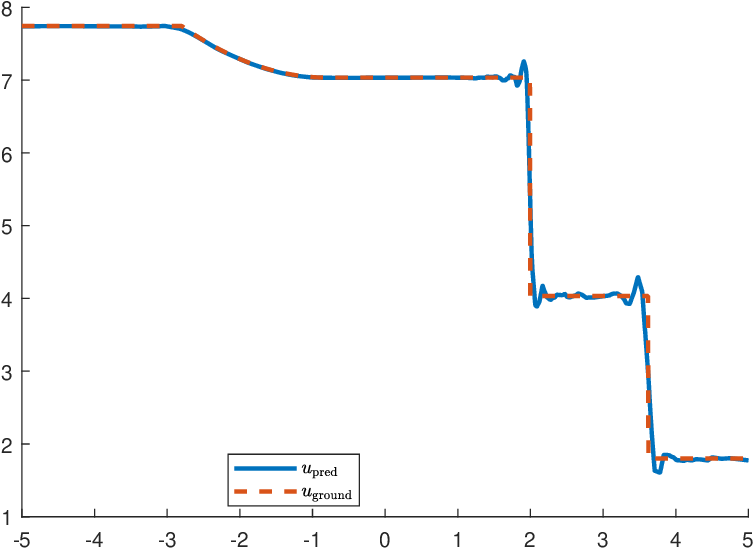}
}
\hspace{2mm}
\subfloat[prediction of Shift DON]{
	\includegraphics[width=0.25\textwidth]{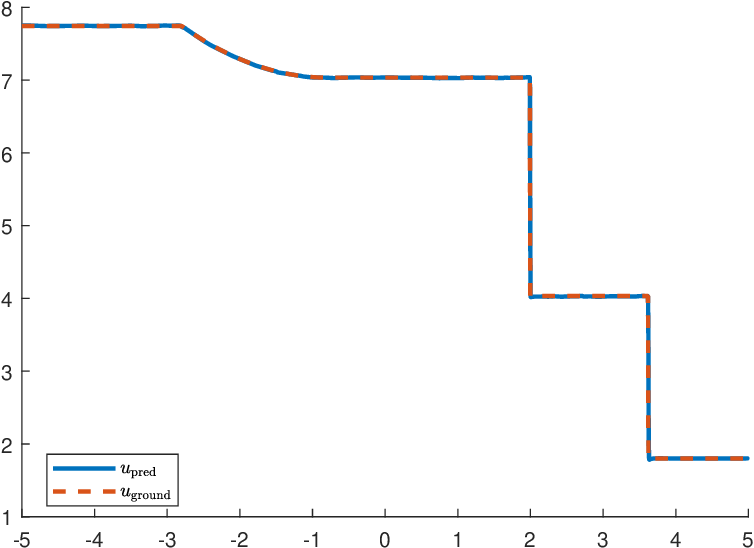}
}
\hspace{2mm}
\subfloat[prediction of R-adaptive DON]{
	\includegraphics[width=0.25\textwidth]{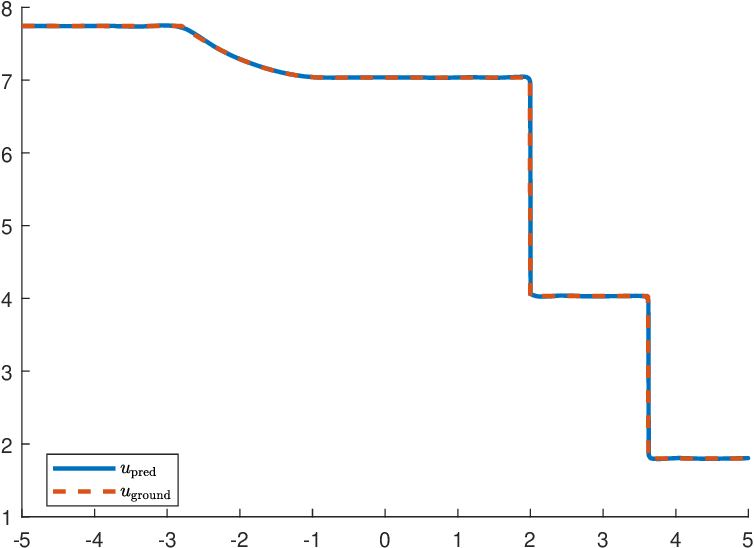}
}

\subfloat[error of vanilla DON]{
	\includegraphics[width=0.25\textwidth]{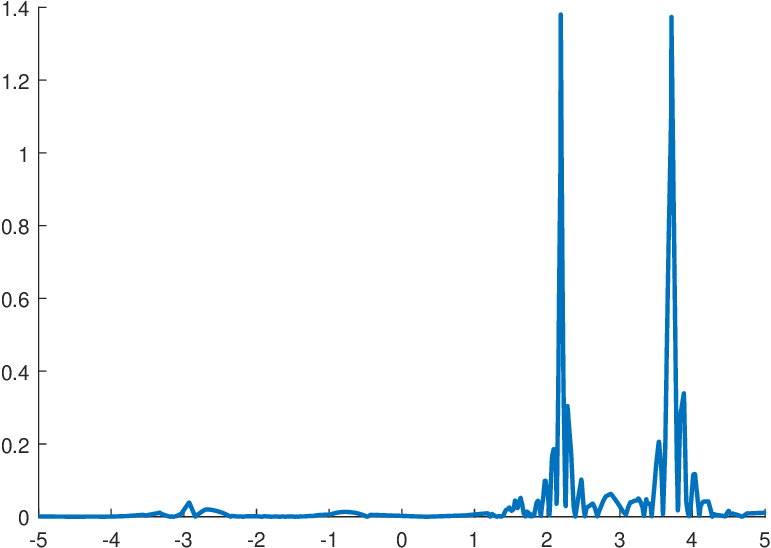}
}
\hspace{2mm}
\subfloat[error of Shift DON]{
	\includegraphics[width=0.25\textwidth]{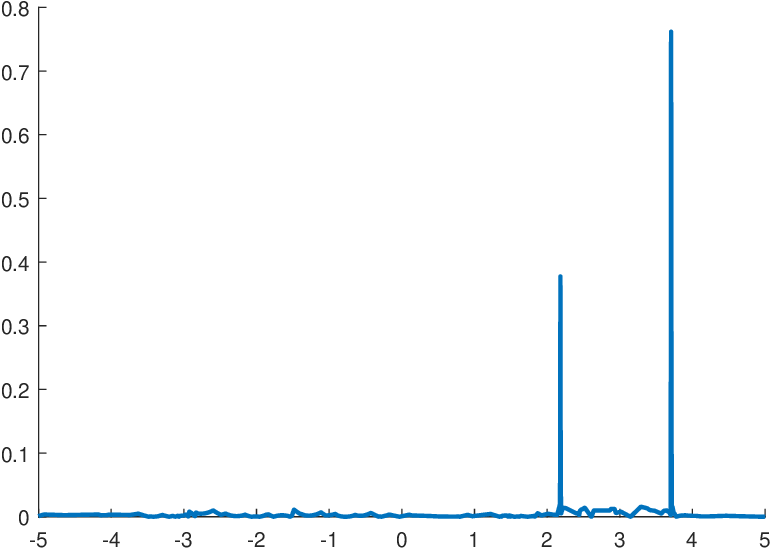}
}
\hspace{2mm}
\subfloat[error of R-adaptive DON]{
	\includegraphics[width=0.25\textwidth]{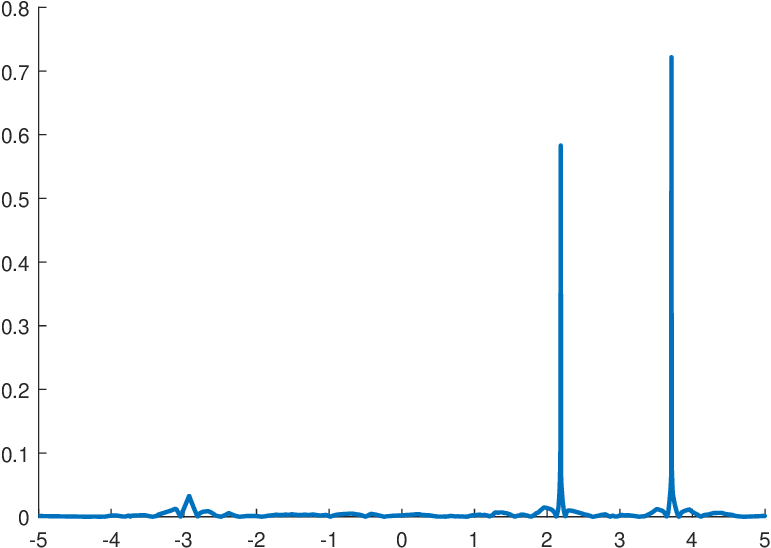}
}
	\caption{Illustration of an example of outputs for Sod shock tube problem.}
	\label{fig:shocktuberesultfigure} 
\end{figure}

In the testing part, for the prediction results of R-adaptive DeepONet, we still obtained them by piecewise linear interpolation to the uniformly distributed grids through the output of two sub-DeepONets. The results are summarized in \cref{tab:shocktubeerrortable} and an example of the output is shown in \cref{fig:shocktuberesultfigure}.
\begin{table}[tbhp]
\footnotesize
\caption{Testing $L^2$ error of different models for Sod shocktube problem.}
\label{tab:shocktubeerrortable}
\centering
\begin{tabular}{cccc}
\toprule[2pt] 
 Model & vanilla DON & Shift-DON & R-adaptive DON  \\
\midrule 
Error & $4.77\times 10^{-4}$ & $2.71\times 10^{-5}$ & $9.24\times 10^{-5}$  \\
\bottomrule[2pt] 
\end{tabular}
\end{table}

From \cref{tab:shocktubeerrortable}, we can see that R-adaptive DeepONet has stronger performance than vanilla DeepONet when approximating the solution of the Sod shock tube problem. And, from \cref{fig:shocktuberesultfigure}, we can see that although the error is larger than that of Shift-DeepONet, R-adaptive DeepONet can catch the discontinuity just as well as Shift-DeepONet.

\section{Conclusion and Discussion}

In this paper, we have proposed a DeepONet learning framework based on the R-adaptive method to address the limitations of vanilla DeepONet representations. Inspired by the introduction of adaptive coordinates in R-adaptive methods, our framework tackles the challenge of representing problems with local singularities by separately learning the adaptive coordinate transform function and the corresponding solution over the computation domain. Additionally, we have derived two solution-dependent weighting strategies in the training process to reduce the final error.

We have established an upper bound on the reconstruction error of DeepONet using error estimation from the piecewise linear interpolation and theoretically demonstrated that our R-adaptive DeepONet framework can reduce this upper bound, indicating its potential for problems with local singularities or discontinuities.

In numerical experiments, we selected several typical partial differential equations with local singularities and used the R-adaptive DeepONet to solve them. We compared its results with those of vanilla DeepONet and Shift-DeepONet. It is shown that R-adaptive DeepONet generally outperforms vanilla DeepONet with smaller approximation errors. Moreover, when compared to Shift-DeepONet, each method had its own advantages depending on the problem. We found that R-adaptive DeepONet can capture the local singularity of the solution with less training data, making it perform better in the space-time Burgers' equation {(see \cref{appendix4}).} 

It is important to note that our proposed R-adaptive DeepONet is a learning framework that does not alter the fundamental structure of DeepONet, making it compatible with other approaches like Shift-DeepONet. Future work will explore potential synergies between these methods to achieve better performance.
\appendix

\section{A brief introduction to R-adaptive method and equidistribution}\label{appendix1}

 In this appendix, we provide a brief introduction to the R-adaptive method and its associated equidistribution principle. More detailed discussions can be found in \cite{MR2722625, MR2538864}.

Suppose that we have a PDE with solution $u(x,t)$, which is posed in a physical domain $\Omega_P\subset \mathbb{R}^d$ with independent spatial variable $x\in \mathbb{R}^d$ for each time $t$. Conceptually, an R-adaptive method generates a moving mesh, continuously mapping a suitable computational space $\Omega_C$ into $\Omega_P$. To achieve this, we assume that a computational coordinate $\xi \in \mathbb{R}^d$ is continuously mapped to the physical coordinate so that $x = x(\xi, t)$.
The basis of the R-adaptive methods is that a fixed set of mesh grids (with fixed connectivity) in $\Omega_C$ is moved by this map to a moving set of grids in $\Omega_P$ where the solution is developing an interesting structure. As a result, a fixed set of basis functions (corresponding to the fixed mesh grids) in $\Omega_C$ is mapped to the adaptive basis functions in $\Omega_P$ for each $t$.      
We write the function in computational coordinate $\xi$ corresponding to $u(x,t)$ as $\tilde{u}(\xi, t) = u(x(\xi, t), t)$.
The structure of the function set $\{\tilde{u}(\xi, t)\}_t$ is much less complicated than $\{u(x, t)\}_t$, allowing us to linearly reconstruct $\{\tilde{u}(\xi, t)\}_t$ with fewer basis functions than $\{u(x, t)\}_t$.

The equidistribution principle plays a fundamental role in the mesh adaptation process. This concept, originating from de Boor \cite{MR0431606}, is a powerful method for identifying a suitable mapping. To implement it, we introduce a (time-dependent) Stieltjes measure $\rho(x,t){\rm d}x$ into the physical domain. The scalar function $\rho(x,t)>0$, known as the mesh density specification function (or monitor function), is designed to be large in regions of $\Omega_P$ where the mesh grids need to be clustered. This function is often defined indirectly via the solution, such that $\rho(x,t) = \rho(x, u(x,t), \nabla u(x,t), \dots, t)$. We do not consider the specific choice of the function $\rho$ here. More detailed discussions can be found in \cite{MR2722625}.

Now introduce an arbitrary non-empty set $K\subset \Omega_C$ in the computational domain, with a corresponding image set $x(K,t)\subset \Omega_P$. The map $x(\cdot, t)$ equidistributes the respective density function $\rho$ if the Stieltjes measure of $K$ and $x(K,t)$ normalized over the measure of their respective domains are the same. This implies that
$$
\frac{\int_K{\rm d}\xi}{\int_{\Omega_C}{\rm d}\xi} = 
\frac{\int_{x(K,t)}\rho(x,t){\rm d}x}{\int_{\Omega_P}\rho(x,t){\rm d}x}.
$$
It follows from a change of variables that
$$
\frac{\int_K{\rm d}\xi}{\int_{\Omega_C}{\rm d}\xi} = 
\frac{\int_K \rho(x(\xi, t),t)|J(\xi, t)|{\rm d}\xi}{\int_{\Omega_P}\rho(x, t){\rm d}x},
$$
where 
$$
J(\xi, t) =  {\rm det}\left(\frac{\partial x(\xi, t)}{\partial \xi} \right).
$$
As the set $K$ is arbitrary, the map $x(\xi, t)$ must obey the identity
\begin{equation}
\label{eq:equidistribution}
\rho(x(\xi,t), t) |J(\xi, t)| =\frac{\int_{\Omega_P}\rho(x, t){\rm d}x}{\int_{\Omega_C}{\rm d}\xi} =: \sigma(t).
\end{equation}
We shall refer to \cref{eq:equidistribution} as the equidistribution equation, and it must always be satisfied by the map $x(\xi, t)$.
By solving the mesh equation \cref{eq:equidistribution} and the original problem simultaneously, we can obtain the adaptive mesh and the corresponding adaptive solution.

\section{Proof of \cref{thm:linearanalysis}}
\label{appendix2}

Recall that with initial data $\bar{u}_{\zeta} = \bar{u}_0(\cdot -\zeta)$, the solution at $t = T$ can be written as 
$$
\mathcal{G}(\bar{u}_{\zeta})(x) = \bar{u}(x - aT - \zeta) = \mathbbm{1}_{[-\frac{\pi}{2}, \frac{\pi}{2}]}(x - aT -\zeta).
$$
Given $\delta>0$, let 
$$
\begin{aligned}
\mathcal{G}_{\delta}(\bar{u}_{\zeta})(x) =
& \frac{1}{\delta}\sigma(x +\frac{\pi}{2} + \frac{\delta}{2} -aT-\zeta) - \frac{1}{\delta}\sigma(x +\frac{\pi}{2} - \frac{\delta}{2}-aT-\zeta)\\
&- \frac{1}{\delta}\sigma(x -\frac{\pi}{2} + \frac{\delta}{2}-aT-\zeta) + \frac{1}{\delta}\sigma(x -\frac{\pi}{2} - \frac{\delta}{2}-aT-\zeta),
\end{aligned}
$$
where $\sigma$ is the rectified linear unit (ReLU). We have that $\mathcal{G}_{\delta} \to \mathcal{G}$ as $\delta \to 0$, or
\eq{\label{appe1}
\left\Vert \mathcal{G}_{\delta}(\bar{u}) - \mathcal{G}(\bar{u})\right\Vert_2^2 
= 4\int_0^{\frac{\delta}{2}} \left( \frac{x}{\delta} \right)^2 \text{d} x
= \delta/6\quad \forall \bar{u}\sim \mu.
}
Since $\delta$ is arbitrary, we can try to approximate $\mathcal{G}_{\delta}$ instead of $\mathcal{G}$. We divide the proof into the following four steps:

\noindent
\textbf{Step 1:}
In the first step, we divide the object operator $\mathcal{G}_{\delta}$ into two parts.
For each $\bar{u}\sim \mu$, we introduce a coordinate transform $x = x(\xi): [-\pi, \pi] \to [-\pi, \pi]$ to its image function $\mathcal{G}_{\delta}(\bar{u})$, referred to as $u$ for convenience, satisfying the equidistribution relation
$$
\left( \rho x_{\xi} \right)_{\xi} = 0,
$$
for the mesh density function
$$
\rho(x) = \sqrt{1 + (\pi^2 - 2\pi\delta) u_x^2},
$$
and the boundary conditions $x(-\pi) = -\pi, x(\pi) = \pi$. So we can get the object function in the transformed variable $\tilde{u}(\xi) = u(x(\xi)).$
For example, when $\zeta = -aT$, i.e. 
$$
\begin{aligned}
u(x) =
& \frac{1}{\delta}\sigma(x +\frac{\pi}{2} + \frac{\delta}{2}) - \frac{1}{\delta}\sigma(x +\frac{\pi}{2} - \frac{\delta}{2})\\
&- \frac{1}{\delta}\sigma(x -\frac{\pi}{2} + \frac{\delta}{2}) + \frac{1}{\delta}\sigma(x -\frac{\pi}{2} - \frac{\delta}{2}),
\end{aligned}
$$
we have that
$$
\begin{aligned}
x(\xi) = -\pi & + (2- \frac{2\delta}{\pi})(\xi + \pi) 
+ ( \frac{4\delta}{\pi} - 2)\sigma(\xi +\frac{3\pi}{4}) + ( 2- \frac{4\delta}{\pi} )\sigma(\xi +\frac{\pi}{4})\\
&+ ( \frac{4\delta}{\pi} - 2)\sigma(\xi -\frac{\pi}{4}) + ( 2- \frac{4\delta}{\pi} )\sigma(\xi -\frac{3\pi}{4}),
\end{aligned}
$$
and
$$
\begin{aligned}
\tilde{u}(\xi) =
& \frac{2}{\pi}\sigma(\xi +\frac{3\pi}{4}) -\frac{2}{\pi}\sigma(\xi +\frac{\pi}{4})\\
&- \frac{2}{\pi}\sigma(\xi -\frac{\pi}{4}) + \frac{2}{\pi}\sigma(\xi -\frac{3\pi}{4}).
\end{aligned}
$$
Note that for each $\bar{u}$, there is a unique $\tilde{u}(\xi)$ and a strictly increasing
 $x(\xi)$ corresponding to it. Let us call these two mappings $\tilde{\mathcal{G}}$ and $\mathcal{T}$, respectively, as
$$
\tilde{\mathcal{G}}: \bar{u} \mapsto \tilde{u}(\xi),\quad 
\mathcal{T}: \bar{u} \mapsto x(\xi).
$$
So the objective operator $\mathcal{G}_{\delta}$ is divided into two parts,
$$
\mathcal{G}_{\delta}(\bar{u}) = \tilde{\mathcal{G}}(\bar{u}) \circ (\mathcal{T}(\bar{u}))^{-1}.
$$

\noindent
\textbf{Step 2:}
In this step we will show that we can approximate the operator $\tilde{\mathcal{G}}$ by a DeepONet with $n+1$ trunk and branch net basis functions, with error
$$
\|\tilde{\mathcal{G}}(\bar{u}) - \tilde{\mathcal{G}}_{\theta_G}(\bar{u}) \|_{L^2{[-\pi, \pi]}} \ \lesssim\ n^{-3/2}.
$$

We use {the} finite element interpolation to prove this. Let $\{\xi_i = -\pi + i \frac{2\pi}{n}\}_{i=0}^n$ be the uniform grid nodes on $[-\pi, \pi]$, and $\{\phi_i\}$ be the corresponding piecewise-linear basis functions. For a given $\tilde{u}(\xi)$, define its finite element interpolation $\tilde{u}_I$ as
$$
\tilde{u}_I = \sum_{i=0}^n \tilde{u}(\xi_i) \phi_i.
$$
Note that $\tilde{u}$ itself is piece-linear, with {four corner points} that are $\pi/2$ apart. So $\tilde{u}_I$ is equal to $\tilde{u}$ in the intervals without corner points. 
Suppose that a corner point is $\xi_j + \eta \in [\xi_j, \xi_{j+1}]$, without loss of generality, we assume that the slope on its left is $0$ and right is $2/\pi$.
The $L^2$ error on $[\xi_j, \xi_{j+1}]$ can be estimated as
$$
\begin{aligned}
\|\tilde{u} - \tilde{u}_I\|_{L^2_{[\xi_j, \xi_{j+1}]}}^2
& = \int_{\xi_j}^{\xi_{j+1}}|\tilde{u}(\xi) - \tilde{u}_I(\xi)|^2 {\rm d}\xi\\
& = \int_0^{\eta}|\frac{n}{\pi^2}(\frac{2\pi}{n} - \eta) \tilde{\xi}|^2 {\rm d}\tilde{\xi}
+ \int_{\eta}^{\frac{2\pi}{n}}|\frac{n}{\pi^2}(\frac{2\pi}{n} - \eta) \tilde{\xi} - \frac{2}{\pi}(\tilde{\xi} - \eta)|^2 {\rm d}\tilde{\xi}\\
& = \frac{1}{3}\frac{n^2}{\pi^4}(\frac{2\pi}{n} - \eta)^2 \eta^3 + \frac{1}{3}\frac{n^2}{\pi^4}(\frac{2\pi}{n} - \eta)^3 \eta^2\\
& = \frac{2}{3\pi^3} n (\frac{2\pi}{n} - \eta)^2 \eta^2 \\
&\leq \frac{2\pi}{3} n^{-3}.
\end{aligned}
$$
So we have
{$$
\|\tilde{u} - \tilde{u}_I\|_{L^2{[-\pi, \pi]}} \leq C n^{-3/2},
$$}
where $C = \sqrt{8\pi/3}$ since there are only four corner points.

Since piecewise linear functions can be represented by ReLU neural networks, there is a neural network $\tau: \mathbb{R} \to \mathbb{R}^{n+1}$, mapping $\xi$ to $(\phi_0(\xi), \dots, \phi_n(\xi))$ precisely. And by the universal approximation property of the neural networks, given arbitrary $\varepsilon >0$,  there exists a neural network $\beta : L^1(\mathbb{T}) \cup L^{\infty}(\mathbb{T}) \to \mathbb{R}^{n+1}$, such that
$$
\mathop{\rm sup}_{\bar{u}\sim \mu} \|\beta(\bar{u}) - (\tilde{u}(\xi_0), \dots, \tilde{u}(\xi_n)) \|_{{l^{\infty}}} < \varepsilon.
$$
Then,
$$
\begin{aligned}
\|\tilde{\mathcal{G}}(\bar{u}) - \beta(\bar{u}) \cdot \tau(\cdot) \|_{L^2{[-\pi, \pi]}}
&\leq \|\tilde{u} - \tilde{u}_I\|_2 + \|\tilde{u}_I - \beta(\bar{u}) \cdot \tau(\cdot) \|_2\\
&\lesssim n^{-3/2} + \frac{\varepsilon}{n} \\
&\lesssim n^{-3/2}.
\end{aligned}
$$
{Let $\tilde{\mathcal{G}}_{\theta_G}$ be a DeepONet with trunk net $\tau$ and branch net $\beta$, then we have
\eq{\label{appe3}
\|\tilde{\mathcal{G}}(\bar{u}) -  \tilde{\mathcal{G}}_{\theta_G}(\bar{u}) \|_{L^2{[-\pi, \pi]}} \ \lesssim\ n^{-3/2}.
}}
\noindent
\textbf{Step 3:}
Following the step above, we will show that we can approximate the operator $\mathcal{T}$ by a DeepONet with $n+1$ trunk and branch net basis functions, with error
$$
\|(\mathcal{T}(\bar{u}))^{-1} - (\mathcal{T}_{\theta_T}(\bar{u}))^{-1} \|_{L^2{[-\pi, \pi]}} \ \lesssim\ n^{-3/2}.
$$

Similarly to Step 2, we use {the} finite element interpolation to prove this. For a given $x(\xi)$, define its finite element interpolation $x_I$ as
$$
x_I(\xi) = \sum_{i=0}^n x(\xi_i) \phi_i(\xi).
$$
Since $x_I(\xi)$ is strictly increasing, we denote its inverse as $\xi_I(x)$. $\xi_I(x)$ is linear in each interval $[x_i, x_{i+1}]$, where $x_i = x(\xi_i)$, and $\xi_I(x_i) = \xi(x_i) = \xi_i$ for each $i$, where $\xi(x)$ is the inverse of $x(\xi)$. So $\xi_I(x)$ is equal to $\xi(x)$ in those $[x_i, x_{i+1}]$ without corner points. Suppose a corner point of $x(\xi)$ is $\xi_j + \eta \in [\xi_j, \xi_{j+1}]$, without loss of generality, we assume that the slope on its left is $2-2\delta/\pi$ and right is $2\delta/\pi$.

For simplicity, let $k_1 = 2-2\delta/\pi, k_2 = 2\delta/\pi$. We can calculate that the corner point is $(\xi_j + \eta, x_j + k_1 \eta)$, and $x_{j+1} = x_j + k_1 \eta + k_2 (2\pi/p - \eta)$. So we have
$$
\begin{aligned}
x(\xi) &= \left\{
\begin{aligned}
&k_1 (\xi- \xi_j) + x_j, &&\ \xi\in [\xi_j,\ \xi_j + \eta]{,} \\
&k_2 (\xi- \xi_j - \eta) + x_j + k_1 \eta, &&\ \xi\in [\xi_j + \eta,\ \xi_{j+1}]{,}
\end{aligned}
\right.\\
x_I(\xi) &= k_3 (\xi- \xi_j) + x_j, \quad \xi\in [\xi_j,\ \xi_{j+1}],
\end{aligned}
$$
where $k_3 = \frac{n}{2\pi}(k_1 \eta + k_2(\frac{2\pi}{n}- \eta))$. Correspondingly,
$$
\begin{aligned}
\xi(x) &= \left\{
\begin{aligned}
&\frac{1}{k_1} (x - x_j) + \xi_j, &&\ x\in [x_j,\ x_j + k_1\eta]{,} \\
&\frac{1}{k_2} (x - x_j - k_1\eta) + \xi_j + \eta, &&\ x\in [x_j + k_1\eta,\ x_{j+1}]{,} 
\end{aligned}
\right.\\
\xi_I(x) &= \frac{1}{k_3} (x - x_j) + \xi_j, \quad x\in [x_j,\ x_{j+1}],
\end{aligned}
$$
Then the $L^2$ error on $[x_j, x_{j+1}]$ with respect to $x$ can be estimated as
$$
\begin{aligned}
\|\xi(x) &- \xi_I(x)\|_{L^2_{[x_j, x_{j+1}]}}^2
= \int_{x_j}^{x_{j+1}}|\xi(x) - \xi_I(x)|^2 {\rm d}x\\
&= \int_{x_j}^{x_j + k_1\eta} |(\frac{1}{k_1} - \frac{1}{k_3})(x - x_j)|^2 {\rm d}x + \int_{x_j + k_1\eta}^{x_{j+1}}|(\frac{1}{k_2} - \frac{1}{k_3})(x - x_j - k_1 \eta)|^2 {\rm d}x \\
&= \int_{0}^{k_1\eta} |(\frac{1}{k_1} - \frac{1}{k_3})\tilde{x}|^2 {\rm d}\tilde{x} 
+ \int_{0}^{k_2 (\frac{2\pi}{n} - \eta)} |(\frac{1}{k_2} - \frac{1}{k_3})\tilde{x}|^2 {\rm d}\tilde{x}\\
&= \frac{1}{3}\left( (\frac{1}{k_1} - \frac{1}{k_3})^2 (k_1\eta)^3 + (\frac{1}{k_2} - \frac{1}{k_3})^2 (k_2 (2\frac{2\pi}{n} - \eta))^3\right)\\
&= \frac{1}{3}\frac{(k_1-k_2)^2 \eta^2 (\frac{2\pi}{n}-\eta)^2}{k_1\eta + k_2(\frac{2\pi}{n}-\eta)}\leq  \frac{1}{3} \frac{(k_1-k_2)^2}{k_1}\eta (\frac{2\pi}{n}-\eta)^2\\
&\leq  \frac{1}{3} \frac{(k_1-k_2)^2}{k_1} \frac{4}{27} (\frac{2\pi}{n})^3.
\end{aligned}
$$
Since $\frac{(k_1-k_2)^2}{k_1} \to 2\pi$ as $\delta \to 0$, we have the result
$$
\|\xi(x) - \xi_I(x)\|_{L^2_{[x_j, x_{j+1}]}} \ \lesssim\ n^{-3/2}.
$$
Similar to Step 2, we can build a DeepONet $\mathcal{T}_{\theta_T}$, such that
\eq{\label{appe2}
\|(\mathcal{T}(\bar{u}))^{-1} - (\mathcal{T}_{\theta_T}(\bar{u}))^{-1} \|_{L^2_{[-\pi, \pi]}} \ \lesssim\ n^{-3/2}.
}

\noindent
\textbf{Step 4:}
{In the last step, we will combine the previous conclusions to arrive at the final result. It is obvious that
$$
\begin{aligned}
\|\mathcal{G}(\bar{u}) &- \tilde{\mathcal{G}}_{\theta_G} (\bar{u}) \circ (\mathcal{T}_{\theta_T}(\bar{u}))^{-1} \|_{L^2}
\leq  \|\mathcal{G}(\bar{u}) - \mathcal{G}_{\delta}(\bar{u}) \|_{L^2} \\
&+ \|\mathcal{G}_{\delta}(\bar{u}) - \tilde{\mathcal{G}}_{\theta_G} (\bar{u}) \circ (\mathcal{T}_{\theta_T}(\bar{u}))^{-1} \|_{L^2} \\
= & \|\mathcal{G}(\bar{u}) - \mathcal{G}_{\delta}(\bar{u}) \|_{L^2} + \|\tilde{\mathcal{G}}(\bar{u}) \circ (\mathcal{T}(\bar{u}))^{-1} - \tilde{\mathcal{G}}_{\theta_G} (\bar{u}) \circ (\mathcal{T}_{\theta_T}(\bar{u}))^{-1} \|_{L^2}\\
\leq & \|\mathcal{G}(\bar{u}) - \mathcal{G}_{\delta}(\bar{u}) \|_{L^2} +  \|\tilde{\mathcal{G}}(\bar{u}) \circ (\mathcal{T}(\bar{u}))^{-1} - \tilde{\mathcal{G}}(\bar{u}) \circ (\mathcal{T}_{\theta_T}(\bar{u}))^{-1} \|_{L^2} \\
&+ \|\tilde{\mathcal{G}}(\bar{u}) \circ (\mathcal{T}_{\theta_T}(\bar{u}))^{-1} - \tilde{\mathcal{G}}_{\theta_G} (\bar{u}) \circ (\mathcal{T}_{\theta_T}(\bar{u}))^{-1} \|_{L^2}=:\mathrm{I}_1 +\mathrm{I}_2 +\mathrm{I}_3.
\end{aligned}
$$
From (\ref{appe1}), it follows that $\mathrm{I}_1=\delta/6$. For $\mathrm{I}_2$,  from (\ref{appe2}), it follows that
$$
\mathrm{I}_2 \leq {\rm Lip}(\tilde{\mathcal{G}}(\bar{u})(\xi))\|(\mathcal{T}(\bar{u}))^{-1}(x) - (\mathcal{T}_{\theta_T}(\bar{u}))^{-1}(x) \|_{L^2}
 \ \lesssim\ \frac{2}{\pi} n^{-3/2}.
$$
Further, for $\mathrm{I}_3$, from (\ref{appe3}), it follows that
$$
\begin{aligned}
\mathrm{I}_3
& = \int_{\mathbb{T}}|\tilde{\mathcal{G}}(\bar{u}) \circ (\mathcal{T}_{\theta_T}(\bar{u}))^{-1}(x) - \tilde{\mathcal{G}}_{\theta_G} (\bar{u}) \circ (\mathcal{T}_{\theta_T}(\bar{u}))^{-1}(x) |^2 {\rm d} x\\
& = \int_{\mathbb{T}}|\tilde{\mathcal{G}}(\bar{u})(\xi) - \tilde{\mathcal{G}}_{\theta_G} (\bar{u}) (\xi) |^2 (\mathcal{T}_{\theta_T}(\bar{u}))_{\xi} {\rm d} \xi \\
& \leq (2 - \frac{2}{\pi}\delta) \|\tilde{\mathcal{G}}(\bar{u}) - \tilde{\mathcal{G}}_{\theta_G}(\bar{u}) \|_{L^2}^2 \\
& \ \lesssim\ (2 - \frac{2}{\pi}\delta) n^{-3}.
\end{aligned}
$$
Hence,
$$
\|\mathcal{G}(\bar{u}) - \tilde{\mathcal{G}}_{\theta_G} (\bar{u}) \circ (\mathcal{T}_{\theta_T}(\bar{u}))^{-1} \|_{L^2}
 \ \lesssim\ \delta/6 + \sqrt{(2 - \frac{2}{\pi}\delta)} n^{-3/2} + {\frac{2}{\pi} n^{-3/2}.}
$$
Since $\delta$ is arbitrary, we have the final result
$$
\|\mathcal{G}(\bar{u}) - \tilde{\mathcal{G}}_{\theta_G} (\bar{u}) \circ (\mathcal{T}_{\theta_T}(\bar{u}))^{-1} \|_{L^2}
\lesssim n^{-3/2} \quad \forall \bar{u}\sim \mu.
$$}
This implies the claim of Theorem.


\section {Algorithm of R-adaptive DeepONet and definitions of errors}\label{appendix0}

In this appendix, we summarize the proposed algorithm using in \cref{sec:5}:
\begin{algorithm}[htbp]
	\caption{Framework of R-Adaptive DeepONet}
	\label{alg:algorithm1}
\begin{algorithmic}[1]
\STATE\textbf{Require}: Dataset $\left\{a^k, u^k \right\}_{k=1}^{N_1}$. Here, $u^k$ can be represented in any discrete form, as long as sufficient information about $u^k$ is available.
\STATE\textbf{Data processing}: Apply the equidistributing principle to each $u^k$ to get processed data sets $\{\tilde{y}_j^k = y^k(\xi_j)\}_{j=1}^{N_2}$ and  $\{\tilde{u}_j^k = u^k(\tilde{y}_j^k)\}_{j=1}^{N_2}$. Here $\{\xi_j\}_{j=1}^{N_2}$ is a set of uniformly distributed points on the computational domain.
\STATE\textbf{Training}: Train the adaptive coordinate DON $\mathcal{T}_{\theta_T}: a \mapsto y(\xi)$ using dataset $\{a^k, \{\tilde{y}_j^k\}_{j=1}^{N_2} \}_{k=1}^{N_1}$ through minimizing the loss function $\mathcal{L}_{\mathcal{T}}$ defined in \cref{eq:lossT}.
Train the adaptive solution DON $\tilde{\mathcal{G}}_{\theta_G}: a\mapsto\tilde{u}(\xi)$ using dataset $\{a^k, \{\tilde{u}_j^k\}_{j=1}^{N_2} \}_{k=1}^{N_1}$ through minimizing the loss function $\mathcal{L}_{\tilde{\mathcal{G}}}$ defined in \cref{eq:lossG}.
\STATE\textbf{Predicting}: 
Given a new parameter $a\in \mathcal{X}$ and points $y\in D_{\mathcal{Y}}$.
Compute the adaptive mesh $y_j = \mathcal{T}_{\theta_T}(a)(\xi_j)$ and the corresponding solution values $u(y_j) = \tilde{\mathcal{G}}_{\theta_G}(a)(\xi_j)$. 
Based on the obtained discrete points $\{y_j, u(y_j)\}$, determine the solution values $u(y)$ using methods such as the interpolation method.
\end{algorithmic}
\end{algorithm}

{We use the relative $L^2$ norm $E_r$ as the error metric employed throughout all numerical experiments to assess model performance. Specifically, the reported testing errors correspond to the mean of the absolute and relative $L^2$ error of a trained model over all examples in the test data-set, i.e.
$$
E_r := \frac{1}{N_1} \sum_{k=1}^{N_1} \frac{\|\mathcal{G}_{\theta}(a_k) - \mathcal{G}(a_k)\|_2}{\|\mathcal{G}(a_k)\|_2},
$$
where $\mathcal{G}$ is the target operator and $\mathcal{G}_{\theta}$ is the model to be evaluated, $N_1$ is the number of examples in the test data-set. The $L^2$ norm is discretely computed through
$$
\|\mathcal{G}(a) \|_2 = \sqrt{\frac{1}{N_2} \sum_{j=1}^{N_2} |\mathcal{G}(a)(y_j) - \mathcal{G}_{\theta}(a)(y_j)|^2},
$$
where $\{y_j\}_{j=1}^{N_2}$ is a set of equi-distributed points in the domain of $\mathcal{G}(a)$.}

\section{Experiment of Spatial-temporal Burgers' equation}\label{appendix4}

To verify that our proposed framework is suitable for higher dimensional cases, we extend the previous experiment to learn the operator mapping initial conditions $u(x,0)$ to the spatial-temporal solution $u(x,t)$ of the viscous Burgers' Equation \cref{eq:burgersequation} with fixed viscosity $\nu = 10^{-4}$.

\begin{figure}[tbhp] 
\centering
\subfloat[ground solution]{
	\includegraphics[width=0.25\textwidth]{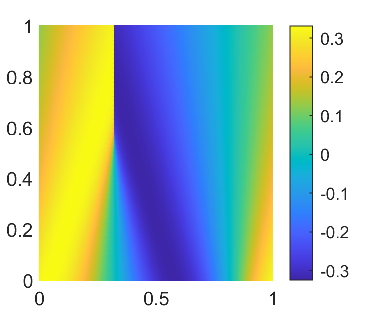}
}
\hspace{2mm}
\subfloat[adaptive coordinate]{
	\includegraphics[width=0.25\textwidth]{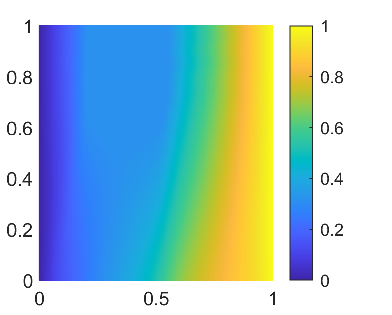}
}
\hspace{2mm}
\subfloat[adaptive solution]{
	\includegraphics[width=0.25\textwidth]{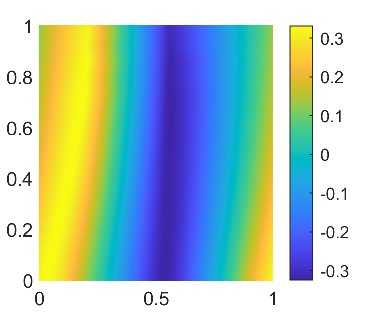}
}

\hspace{2mm}
\subfloat[weight $w_{\mathcal{T}}$]{
	\includegraphics[width=0.25\textwidth]{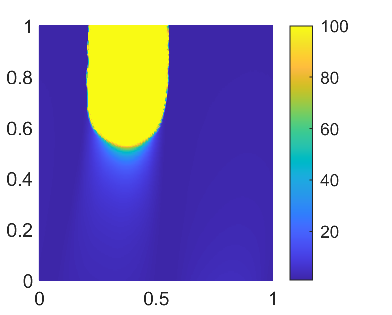}
}
\hspace{2mm}
\subfloat[weight $w_{\tilde{\mathcal{G}}}$]{
	\includegraphics[width=0.25\textwidth]{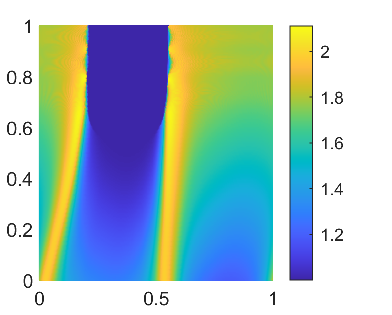}
}

	\caption{An example of the processed data for Spatial-temporal Burgers' equation.}
	\label{fig:burgers2ddata} 
\end{figure}

The training and testing data are generated as that of previous experiment. The training output are generated by sampling the underlying exact solution on a spatial-temporal mesh of resolution $257\times 201$. Considering that such functions generally have no singularity in the temporal direction, we only perform equal distribution  in the spatial direction during data preprocessing, and the density function used is still \begin{equation} 
\rho(x,t) = \sqrt{1+\left|\frac{\partial u}{\partial x}\right|^2}.
\end{equation}
An example of the processed data is presented in \cref{fig:burgers2ddata}. It can be seen that the adaptive coordinates and adaptive solutions obtained by preprocessing do not have the obvious singularity {compared to} the ground solution. This indicates that we can approximate them well with the sub-DeepONet separately.

For such a $2$-d problem, we use a model with a larger size. For the vanilla DeepONet, the two sub-DeepONets of R-adaptive DeepONet, and the main body of Shift-DeepONet, we use the same structure: both the branch net and trunk net have 5 layers with 256 neurons each. For Shift-DeepONet's scale net and shift net, we also used a structure of 5 layers with 256 neurons per layer to keep the number of parameters in the model as close to the same order as possible. The training process is set up similarly to the previous subsection. 

\begin{figure}[tbhp] 
\centering
\subfloat[prediction of vanilla DON]{
	\includegraphics[width=0.25\textwidth]{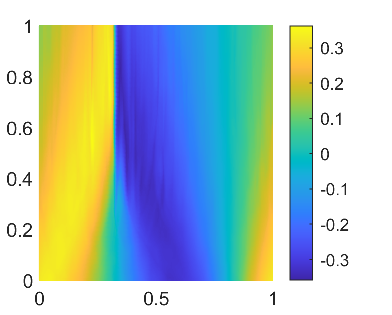}
}
\hspace{2mm}
\subfloat[prediction of Shift DON]{
	\includegraphics[width=0.25\textwidth]{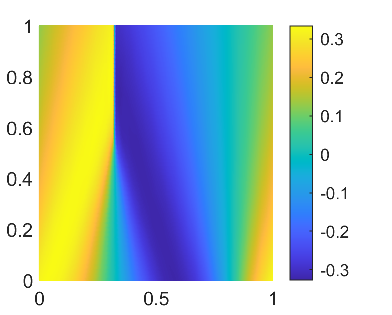}
}
\hspace{2mm}
\subfloat[prediction of R-adaptive DON]{
	\includegraphics[width=0.25\textwidth]{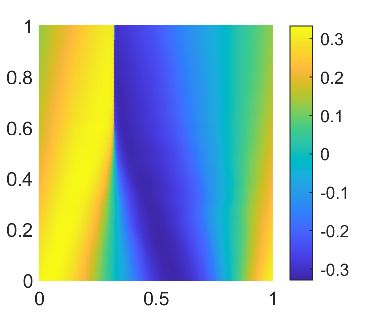}
}

\subfloat[error of vanilla DON]{
	\includegraphics[width=0.25\textwidth]{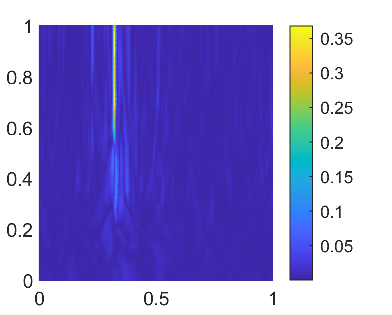}
}
\hspace{2mm}
\subfloat[error of Shift DON]{
	\includegraphics[width=0.25\textwidth]{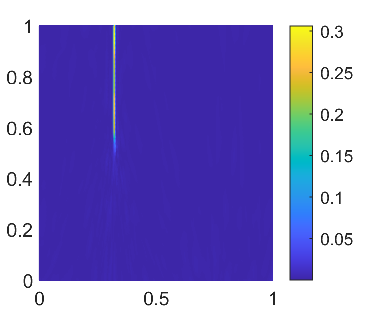}
}
\hspace{2mm}
\subfloat[error of R-adaptive DON]{
	\includegraphics[width=0.25\textwidth]{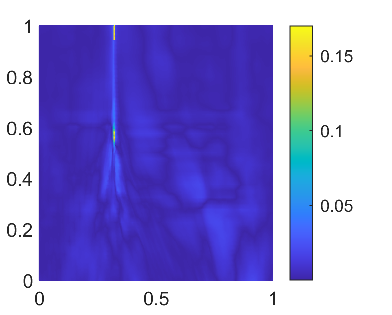}
}
	\caption{Illustration of an example of outputs for spatial-temporal Burgers' equation.}
	\label{fig:burgers2dresult} 
\end{figure}

In the test part, we use vanilla DeepONet and Shift-DeepONet to directly predict the function values on uniformally distributed spatial-temporal grids with size $257\times201$, and compare them with the ground solution. For R-adaptive DeepONet, we use the trained adaptive solution DeepONet and adaptive coordinate DeepONet to output the coordinates of adaptive distributed grids of size $65\times41$ and their corresponding solution values. {And then we utilize the bilinear interpolation to calculate function values on $257\times201$ uniformally distributed grids, which are used to compare with the ground solution. }The testing errors are shown in \cref{tab:burgers2derrortable} and an output example is shown in \cref{fig:burgers2dresult}. {It can be seen that in this experiment, R-adaptive DeepONet successfully reduced the approximation error, and even performed better than the Shift-DeepONet. This may be due to the fact that the insufficient resolution of the training data  is not enough for Shift-DeepONet to grasp local singularity accurately. Comparatively, R-adaptive DeepONet, through data preprocessing, can capture the local singularity of the solution with lower-resolution data. This observation demonstrates another advantage of R-adaptive DeepONet: it can reduce the size of training data through adaptive data sampling. }

\begin{table}[tbhp]
\footnotesize
\caption{Testing $L^2$ error of different models for spatial-temporal Burgers' equation.}
\label{tab:burgers2derrortable}
\centering
\begin{tabular}{ccc}
\toprule[2pt] 
 Model & $E_a$ & $E_r$  \\
\midrule 
vanilla DeepONet & $7.61\times 10^{-4}$ & $1.58\times 10^{-2}$  \\
Shift-DeepONet & $5.04\times 10^{-4}$ & $2.86\times 10^{-2}$ \\
R-adaptive DON & $2.74\times 10^{-4}$ & $3.34\times 10^{-3}$ \\
\bottomrule[2pt] 
\end{tabular}
\end{table}

\bibliographystyle{siamplain}
\bibliography{refer}

\end{document}